\documentclass[12pt]{amsart}
\usepackage{amssymb}
\usepackage{amscd}
\usepackage{graphicx}
\usepackage{pstricks}
\usepackage{pst-plot}

\newtheorem{thm}{Theorem}[section]
\newtheorem{lem}[thm]{Lemma}
\newtheorem{cor}[thm]{Corollary}

\newtheorem{prop}[thm]{Proposition}

\theoremstyle{remark}

\newtheorem{rem}[thm]{Remark}

\theoremstyle{definition}
\newtheorem{defn}[thm]{Definition}

\numberwithin{equation}{section}

\allowdisplaybreaks[4]

\DeclareMathOperator{\Ext}{Ext}
\DeclareMathOperator{\Pic}{Pic}

\DeclareMathOperator{\im}{Im}
\DeclareMathOperator{\ks}{ks}
\DeclareMathOperator{\Aut}{Aut}
\DeclareMathOperator{\Span}{Span}

\begin{document}

\vfuzz0.5pc
\hfuzz0.5pc 

\newcommand{\claimref}[1]{Claim \ref{#1}}
\newcommand{\thmref}[1]{Theorem \ref{#1}}
\newcommand{\propref}[1]{Proposition \ref{#1}}
\newcommand{\lemref}[1]{Lemma \ref{#1}}
\newcommand{\coref}[1]{Corollary \ref{#1}}
\newcommand{\remref}[1]{Remark \ref{#1}}
\newcommand{\conjref}[1]{Conjecture \ref{#1}}
\newcommand{\questionref}[1]{Question \ref{#1}}
\newcommand{\defnref}[1]{Definition \ref{#1}}
\newcommand{\secref}[1]{Sec. \ref{#1}}
\newcommand{\ssecref}[1]{\ref{#1}}
\newcommand{\sssecref}[1]{\ref{#1}}
\newcommand{\figref}[1]{Figure \ref{#1}}
\newcommand{\appendixref}[1]{Appendix \ref{#1}}

\def \mapright#1{\smash{\mathop{\longrightarrow}\limits^{#1}}}
\def \mapleft#1{\smash{\mathop{\longleftarrow}\limits^{#1}}}
\def \mapdown#1{\Big\downarrow\rlap{$\vcenter{\hbox{$\scriptstyle#1$}}$}}
\def \smapdown#1{\downarrow\rlap{$\vcenter{\hbox{$\scriptstyle#1$}}$}}
\def \A{{\mathbb A}}
\def \I{{\mathcal I}}
\def \J{{\mathcal J}}
\def \CO{{\mathcal O}}
\def \C{{\mathcal C}}
\def \BC{{\mathbb C}}
\def \BQ{{\mathbb Q}}
\def \m{{\mathcal M}}
\def \H{{\mathcal H}}
\def \S{{\mathcal S}}
\def \Z{{\mathcal Z}}
\def \BZ{{\mathbb Z}}
\def \W{{\mathcal W}}
\def \Y{{\mathcal Y}}
\def \T{{\mathcal T}}
\def \P{{\mathbb P}}
\def \CP{{\mathcal P}}
\def \G{{\mathbb G}}
\def \F{{\mathbb F}}
\def \BR{{\mathbb R}}
\def \D{{\mathcal D}}
\def \L{{\mathcal L}}
\def \f{{\mathcal F}}
\def \X{{\mathcal X}}

\def \closure#1{\overline{#1}}
\def \EQ{\Leftrightarrow}
\def \imply{\Rightarrow}
\def \isom{\cong}
\def \embed{\hookrightarrow}
\def \tensor{\mathop{\otimes}}
\def \wt#1{{\widetilde{#1}}}


\title{A Simple Proof that Rational Curves on K3 are Nodal}

\author{Xi Chen}

\address{Department of Mathematics\\
South Hall, Room 6607\\
University of California\\
Santa Barbara, CA 93106}
\email{xichen@math.ucsb.edu}
\date{July 15, 2001}
\thanks{Research supported in part by a Ky Fan Postdoctoral Fellowship}
\maketitle

\section{Introduction and Statement of Results}

The purpose of this paper is to give a simple proof of the following
theorem proved in \cite{C2}.

\begin{thm}\label{t1}
All rational curves in the primitive class of a general K3
surface of genus \(g\ge 2\) are nodal.
\end{thm}

Please see \cite{C1} and \cite{C2} for the background of this problem.

We will use a degeneration argument as in \cite{C2}. But instead of
degenerating a general K3 surface to a pair of rational surfaces, we
will specialize it to a K3 surface $S$ with Picard lattice
\begin{equation}\label{s2:e1}
\begin{pmatrix}
-2 & 1\\
1 & 0
\end{pmatrix}.
\end{equation}
The Picard group of $S$ is generated by two effective divisors $C$ and
$F$ with $C^2 = -2$, $F^2 = 0$ and $C\cdot F = 1$. It can be realized
as an elliptic
fibration over $\P^1$ with a unique section $C$, fibers $F$ and
$\lambda = 2$. Here \(\lambda = c_1(\pi_* \omega)\) is the first Chern
class of the Hodge bundle \(\pi_* \omega\) of the fibration \(\pi:
S\to \P^1\) (see \cite{H-M}). It is a standard result that
the number of nodal fibers of an elliptic fibration are given by
\(12\lambda\) \cite[p. 158]{H-M}. So there are exactly \(24\)
rational nodal curves in the linear series \(|F|\) for \(S\) general.

This is the same special K3
surface used by Bryan and Leung in their counting of curves on K3
surfaces \cite{B-L}. It is actually the attempt to understand their
method that leads us to our proof. We will call a K3 surface with
Picard lattice \eqref{s2:e1} a {\it BL K3\/} surface.

A BL K3 surface $S$ lies on the boundary of the moduli space of K3
surfaces of
genus $g$ with $C + gF$ as the corresponding primitive divisor. Every
curve in the linear series $|\CO_S(C + gF)|$ is ``totally reducible'',
i.e., it
consists of the $-2$ curve $C$ and $g$ elliptic ``tails'' attached to
$C$. A curve \(D\in |\CO_S(C + gF)|\) is the image of a stable rational
map only
if $D = C\cup m_1 F_1 \cup m_2 F_2 \cup ...\cup m_{24} F_{24}$, where
$F_1, F_2, ..., F_{24}$ are $24$ rational nodal curves in the pencil
$|F|$ and $\sum_{i=1}^{24} m_i = g$; $D$ is obviously nodal if $m_i
\le 1$ for all $i$. The main problem is, of course,
$m_i$ might be greater than 1, i.e., $D$ might be nonreduced, in which
case we need to show that when $S$ deforms to a general K3 surface
$S'$ of genus $g$ and $D$ correspondingly deforms to a rational curve
$D'\subset S'$, $D'$ is necessarily nodal.

It is worthwhile to mention that although this proof looks quite
different from the one in \cite{C2}, all the basic techniques have
already been developed there. By choosing a ``good'' degeneration as the
one used by Bryan-Leung, we eliminate a substantial amount of
technicality in the previous proof. In addition, this proof also gives
a geometric interpretation of Bryan-Leung's work and makes it
possible to redo their counting in the frame of classical
algebraic geometry, if one chooses so. Indeed, we will recover part of
their counting formula in \appendixref{apdx2}.

We will work exclusively over $\BC$ throughout the paper. We use the
usual topology instead of Zariski topology most of the time. When we
say ``neighborhood'' of a point or a subscheme, we usually mean
analytic neighborhood.

\medskip\noindent{\bf Acknowledgments.}
I came up with the main idea of this paper
during a pleasant visit of UT Austin. I would like to thank
Sean Keel for his invitation and for some very helpful conversations
with him. I am especially grateful to the referee, who provided me a
long and detailed report. His corrections and suggestions help me
improve the paper greatly not only in mathematics but also in exposition.
In particular, all the pictures in the current version were drawn and
supplied to me by the referee, in the hope that they will make the
paper more readable.

\section{Degeneration of K3 surfaces}\label{s2}

Let $X$ be a smooth family of K3 surfaces of genus $g$ over the disk
$\Delta$ whose central fiber $X_0 = S$ is a BL K3. Let $Y\subset X$ be
a flat family of rational curves with $Y_t\subset X_t$ and $Y_0\in
|C+gF|$, where $\Delta$ is parameterized by $t$ and $Y_t$ and $X_t$
are general fibers of $Y$ and $X$ over $t\ne 0$.
Notice that a base change might be needed to ensure the
existence of $Y$. Let $E$ be one of the $24$ rational curves $F_1, F_2,
..., F_{24}$ and $p\in E$ be the node of $E$.
Suppose that $Y_0$ contains $E$ with multiplicity
$m$. It suffices to show that $Y_t$ has $m$ nodes in the neighborhood
of $E$. If $m = 1$, there is nothing to prove; otherwise,
we need to apply the stable reduction to $Y$ by
blowing up $X$ and $Y$ along $E$.

Let $N_{A/B}$ denote the normal bundle of $A\subset B$. Here the
normal bundle is defined as the dual of conormal bundle, i.e.,
\begin{equation}\label{s2:e1.1}
N_{A/B} = {\mathcal H}om(I_A/I_A^2, \CO_A),
\end{equation}
where $I_A$ is the ideal sheaf of $A$ in $B$.

If we blow up $X$ along $E$ (see \figref{F1}),
the exceptional divisor is a ruled
surface over $E$ given by $\P N_{E/X}$. We have the exact sequence
\begin{equation}\label{s2:e2}
0 \xrightarrow{} N_{E/S} \xrightarrow{} N_{E/X} \xrightarrow{}
\left.N_{S/X}\right|_E \xrightarrow{} 0.
\end{equation}

Notice that 
\begin{equation}\label{s2:e2.1}
N_{E/S} = \left.N_{S/X}\right|_E = \CO_E \text{ and }
\Ext(\CO_E, \CO_E) = H^1(\CO_E) = \BC.
\end{equation}
So \eqref{s2:e2} might not split. Actually this is always the case as
long as $X$ is general enough. We claim that

\begin{prop}\label{prop1}
The exact sequence \eqref{s2:e2} does not split provided that
the Kodaira-Spencer class of \(X\) is general.
\end{prop}

\begin{rem}\label{rem1}
Some explanations might be needed on what exactly we mean by a
general Kodaira-Spencer class as stated in the above proposition.
The first order deformations of \(S\) are classified by \(H^1(T_S)\)
and the Kodaira-Spencer map of $X$ is
\begin{equation}\label{rem1:e1}
\ks: T_{\Delta, 0} \isom H^0(N_{S/X})\to H^1(T_S),
\end{equation}
where \(T_{\Delta, 0}\) is the tangent space of \(\Delta\) at the origin and
\(T_S\) is the tangent bundle of $S$.
The versal deformation space of \(S\) as a
complex manifold has dimension \(h^1(T_S) = 20\). However, not every
vector in \(H^1(T_S)\) is the Kodaira-Spencer class of a
projective family \(X\). The algebraic deformations of \(S\)
are actually given by the vectors of \(H^1(T_S)\) lying in a
union of countably many subspaces of codimension 1.
This is a well-known fact. However, we need the following more
precise statement.

\begin{lem}\label{lem0}
Let \(X\) be a smooth family of complex surfaces over \(\Delta\) whose
central fiber \(X_0 = S\) is a surface with trivial canonical bundle.
Let \(Y\subset X\) be a closed subscheme of \(X\) of
codimension 1 which is flat over \(\Delta\) and whose central fiber
\(Y_0 = D\) is an ample divisor on \(S\). Then
the Kodaira-Spencer class \(\ks(\partial/\partial t)\) of \(X\)
lies in the subspace \(V\subset H^1 (T_S)\) consisting of the vectors
which are perpendicular to the first Chern class
\(c_1(D)\in H^1(\Omega_S)\) of the divisor \(D\), i.e.,
\begin{equation}\label{lem0:e1}
\ks(\partial/\partial t) \in V
= \{ v\in H^1(T_S): {<}v, c_1(D){>} = 0\},
\end{equation}
where \(\Omega_S\) is the cotangent sheaf of \(S\) and
the pairing \({<}\cdot, \cdot{>}\) is given by Serre duality
\(H^1(T_S) \times H^1(\Omega_S)\to\BC\).

On the other hand, if we fix a K3 surface \(S\) and an ample divisor
\(D\) on \(S\), then for each \(v\in V\), there
exists a pair \((X, Y)\) such that \(Y\subset X\), \(X_0 = S\), \(Y_0
= D\) and the Kodaira-Spencer class of \(X\) is \(v\).
\end{lem}

We are quite certain that the above lemma is also well known.
But since we are unable to locate a reference for it, we will give a
proof in \appendixref{apdx1}.

Roughly, \lemref{lem0} says that a general deformation of a
surface \(S\) with trivial canonical bundle does not preserve
any ample divisor \(D\) on \(S\). 
As a direct consequence, we see that a general deformation of
an algebraic K3 or abelian surface is no longer algebraic.

Back to our situation and we see that the Kodaira-Spencer class of
\(X\) lies the subspace of \(H^1(T_S)\) perpendicular to
\(c_1(C + gF)\), i.e.,
\begin{equation}\label{s2:e2.5}
\ks(\partial/\partial t)\in V
= \{ v\in H^1(T_S): {<}v, c_1(C + gF){>} = 0\}
\end{equation}
by \lemref{lem0}.
Furthermore, for each \(v\in V\), there exists a family \(X\)
whose Kodaira-Spencer class is given by \(v\).
In \propref{prop1}, by \(\ks(\partial/\partial t)\) being general, we
mean that \(\ks(\partial/\partial t)\) is general in \(V\).
\end{rem}

\begin{proof}[Proof of \propref{prop1}]
The sequence \eqref{s2:e2} splits if and only if the induced map
\begin{equation}\label{s2:e3}
H^0(\left.N_{S/X}\right|_E) \xrightarrow{} H^1(N_{E/S})
\end{equation}
is zero. We have the commutative diagram
\begin{equation}\label{s2:e3.1}
\begin{CD}
0 @>>> T_S|_E @>>> T_X|_E @>>> N_{S/X}|_E @>>> 0 \\
 & & @VVV @VVV @VVV \\
0 @>>> N_{E/S} @>>> N_{E/X} @>>> N_{S/X}|_E @>>> 0
\end{CD}
\end{equation}
and we can naturally identify $H^0(\left.N_{S/X}\right|_E)$ with
$T_{\Delta, 0}$. Therefore, the map \eqref{s2:e3} factors through the
Kodaira-Spencer map $T_{\Delta, 0} \to H^1(T_S)$, the restriction
$H^1(T_S) \to H^1(T_S|_E)$ and the surjection 
\begin{equation}\label{s2:e3.2}
H^1(T_S|_E)\xrightarrow{} H^1(N_{E/S})\xrightarrow{} H^2(T_E) = 0.
\end{equation}
In short, we have
\begin{equation}\label{s2:e4}
H^0(\left.N_{S/X}\right|_E) \isom T_{\Delta, 0} \xrightarrow{\ks} H^1(T_S)
\xrightarrow{} H^1(T_S|_E) \xrightarrow{} H^1(N_{E/S}).
\end{equation}
The last map \(H^1(T_S|_E) \rightarrow H^1(N_{E/S})\) is actually an
isomorphism by the following argument.

By the standard exact sequence
\begin{equation}\label{s2:e5}
0\xrightarrow{} N_{E/S}^\vee \xrightarrow{} \Omega_S|_E
\xrightarrow{} \Omega_E \xrightarrow{} 0,
\end{equation}
we have the exact sequence
\begin{equation}\label{s2:e6}
H^0(N_{E/S}) \xrightarrow{} \Ext(\Omega_E, \CO_E) \xrightarrow{}
H^1(T_S|_E)\xrightarrow{} H^1(N_{E/S}) \xrightarrow{} 0.
\end{equation}
Notice that $H^0(N_{E/S}) = \BC$ classifies the embedded deformations
of $E\subset S$ and $\Ext(\Omega_E, \CO_E) = \BC$
classifies the versal deformations of $E$. To show that $H^0(N_{E/S})$
maps nontrivially to $\Ext(\Omega_E, \CO_E)$, it suffices to show that
as $E$ varies in the pencil $|\CO_S(E)|$, the corresponding
Kodaira-Spencer map to the tangent space of
the versal deformation space of $E$ at the origin is
nontrivial, or equivalently, the map to the versal deformation space
of \(E\) is unramified over the origin.
To see this has to be true, we only need to localize the
problem at the node $p$ of $E$: if the map to the versal deformation space
is ramified over the origin, then \(S\) is locally given by \(xy =
t^\alpha\) at \(p\) for some \(\alpha > 1\); however, this is
impossible since \(S\) is smooth at \(p\).
Therefore, the map
\begin{equation}\label{s2:e6.0}
H^0(N_{E/S}) \to \Ext(\Omega_E, \CO_E)
\end{equation}
is nonzero and hence must be an isomorphism. Thus we conclude that 
\begin{equation}\label{s2:e6.1}
H^1(T_S|_E)\xrightarrow{\sim} H^1(N_{E/S}) = \BC
\end{equation}
is an isomorphism.

We have the exact sequence
\begin{equation}\label{s2:e7}
H^1(T_S(-E)) \xrightarrow{f} H^1(T_S) \xrightarrow{} H^1(T_S|_E) \isom
H^1(N_{E/S}) = \BC.
\end{equation}
Combining \eqref{s2:e4} and \eqref{s2:e7}, we are left to
show that the image of the map \(f: H^1(T_S(-E))\to H^1(T_S)\) does not
contain \(V\subset H^1(T_S)\) as in \eqref{s2:e2.5}. 
We claim that that the image of
\(f: H^1(T_S(-E))\to H^1(T_S)\) is contained in the subspace \(W\)
of \(H^1(T_S)\) perpendicular to \(c_1(E)\), i.e.,
\begin{equation}\label{s2:e7.1}
\im f \subset W = \{ v\in H^1(T_S): {<}v, c_1(E){>} = 0\}.
\end{equation}
By Kodaira-Serre duality, we have the following commutative diagram:
\begin{equation}\label{s2:e7.2}
\begin{CD}
H^1(T_S(-E)) @>\sim>> H^1(\Omega_S(E))^\vee @>\sim>>
H^1(\Omega_S(-E))\\
@VV{f}V @VVV @VV{g}V\\
H^1(T_S) @>\sim>> H^1(\Omega_S)^\vee @>\sim>>
H^1(\Omega_S).
\end{CD}
\end{equation}
So we may identify the map \(f\)
with \(g: H^1(\Omega_S(-E)) \to H^1(\Omega_S)\), which is the same
as
\begin{equation}\label{s2:e7.3}
g: H^{1,1}(\CO_S(-E)) \to H^{1,1}(\CO_S)
\end{equation}
on the Dolbeault cohomologies.
For any \(\psi\in H^{1, 1}(\CO_S(-E))\), we have
\begin{equation}\label{s2:e7.4}
\int_S g(\psi)\wedge c_1(E) = \int_E g(\psi) = 0.
\end{equation}
So \eqref{s2:e7.1} follows.

On the other hand, we have \eqref{s2:e2.5}. It is trivial that \(c_1(C
+ gF)\) and \(c_1(E) = c_1(F)\) are linearly independent in
\(H^1(\Omega_S)\).
So \(W\not\supset V\) and a general Kodaira-Spencer class
\(\ks(\partial/\partial t)\in V\) does not lie in \(W\) and hence
\(\ks(\partial/\partial t)\not\in \im f\).
Therefore, \(\ks(\partial/\partial t)\) maps nontrivially
to \(H^1(T_S|_E)\isom H^1(N_{E/S})\). Consequently, the map
\eqref{s2:e3} is not zero and the sequence \eqref{s2:e2} does not split.
\end{proof}

\begin{defn}\label{defn1}
There are two ruled surfaces $\P W$ over $E$, where $W$ is a rank two
vector bundle over $E$ satisfying the exact sequence
\begin{equation}\label{defn1:e1}
0\xrightarrow{} \CO_E\xrightarrow{} W\xrightarrow{} \CO_E\xrightarrow{} 0.
\end{equation}
The proof of this fact is not hard, it goes exactly as the
classification of the ruled surfaces over an elliptic curve with \(e = 0\)
(see e.g. \cite[V, Theorem 2.15]{Ha}) and we will later give a more
geometrical proof of this fact in \ssecref{s3:1}.
If $W = \CO_E\oplus \CO_E$, we
call $\P W \isom\P^1\times E$ {\it trivial\/}; otherwise if $W$ is
indecomposable, we call $\P W$ {\it twisted\/}. 
\end{defn}
 
Even if the family $X$ we start with is general, we cannot draw the
conclusion that $N_{E/X}$ is indecomposable by \propref{prop1}
yet. The problem is that we have already applied a base change to $X$
to ensure the existence of $Y$. If the degree $\alpha$ of the base
change is greater than 1, the Kodaira-Spencer class of the resulting
family $X$ will vanish; and
if we blow up $X$ along $E$, the exceptional divisor is
simply the trivial ruled surface over $E$. But eventually a twisted ruled
surface over $E$ will show up if we keep blowing up $X$ along $E$. We
explain precisely what we mean by this in the next paragraph.

Let $X^{(1)}$ be the blowup of $X$ along $E_0 = E$ (see \figref{F1}).
The central fiber $X_0^{(1)} = S_0\cup S_1$
consists of the proper transform $S_0$ of
$S$ and a ruled surface $S_1$ over $E_0$. If $S_1$ is twisted, we stop at
$X^{(1)}$. Otherwise, $S_1\isom \P^1\times E_0$ is trivial. Notice
that the total family $X^{(1)}$ acquires a singularity during the
blowup; it has a rational double point $p_1\ne p_0\in F_{p_0} \subset
S_1$ over the node $p_0 = p$ of $E_0$, where $F_{p_0}$ is the fiber of
$S_1\to E_0$ over $p_0$. 

\psset{unit=0.01in,linewidth=1.2pt,dash=3pt 2pt,dotsep=2pt}


\begin{figure}[ht]
\centering
\begin{pspicture}(-170,-460)(170,0)

\rput[lt]{0}(-170,0){%
\begin{pspicture}(-170,-210)(170,0)
\psarc[linestyle=dashed](-176.57,-80){80}{315}{45}
\psarc(-64.43,-80){80}{135}{225}
\psarc[linestyle=dashed](-56.57,-80){80}{315}{45}
\psarc(56.57,-80){80}{135}{225}
\psarc(176.57,-80){80}{135}{225}
\psarc(64.43,-80){80}{315}{45}
\psline(-120,-23.43)(120,-23.43)
\psline(-120,-136.57)(120,-136.57)
\psline{->}(0,-150)(0,-170)
\psline(-120,-185)(120,-185)
\pscircle*(0,-185){5}

\parametricplot{-6.7}{6.7}{0.1 t t t mul 30 sub mul mul
t t mul 100 sub}

\rput{0}(0,-110){{\small \(E_0\)}}
\rput{0}(-10,-70){{\small \(p_0\)}}
\rput{0}(0,-200){{\small \(0\)}}
\rput{0}(130,-185){{\small \(\Delta\)}}
\rput{0}(0,-5){{\small \(S_0\)}}
\rput{0}(150,-5){{\small \(X^{(0)}\)}}

\end{pspicture}
}

\psline[doubleline=true,doublesep=1.5pt,arrowsize=2pt 2]{->}(0,-210)(0,-240)

\rput[lt]{0}(-170,-250){%
\begin{pspicture}(-170,-210)(170,0)

\psarc[linestyle=dashed](-176.57,-110){80}{315}{45}
\psarc(-64.43,-110){80}{135}{225}
\psarc[linestyle=dashed](-56.57,-80){80}{315}{348}
\psarc[linestyle=dashed](-56.57,-80){80}{24}{45}
\psarc(56.57,-80){80}{135}{225}
\psarc(176.57,-50){80}{135}{225}
\psarc(64.43,-50){80}{315}{45}

\psline(-120,-53.43)(120,6.57)
\psline(-120,-166.57)(120,-106.57)

\parametricplot{0}{5.477225575}{0.1 t t t mul 30 sub mul mul
t t mul 100 sub}
\parametricplot{-6.7}{-5.477225575}{0.1 t t t mul 30 sub mul mul
t t mul 100 sub}
\parametricplot[linestyle=dashed]{-5.477225575}{0}{%
0.1 t t t mul 30 sub mul mul
t t mul 100 sub}
\parametricplot[linestyle=dashed]{5.477225575}{6.7}{%
0.1 t t t mul 30 sub mul mul
t t mul 100 sub}

\parametricplot{-6.7}{6.7}{72 0.1 t t t mul 30 sub mul mul add
t t mul 76 sub}

\rput{0}(0,-110){{\small \(E_0\)}}
\rput{0}(-10,-70){{\small \(p_0\)}}
\rput{0}(0,-5){{\small \(S_0\)}}
\rput{0}(0,-200){{\small \(0\)}}
\rput{0}(130,-185){{\small \(\Delta\)}}
\rput{0}(150,25){{\small \(X^{(1)}\)}}
\rput{0}(48,-99){{\small \(S_1\)}}
\rput{0}(48,-69){{\small \(F_{p_0}\)}}
\uput[315](72,-76){{\small \(E_1\)}}
\rput{0}(87,-46){{\small \(p_1\)}}

\psline[linestyle=dashed](9.98,-55.11)(71.36,-34.65)
\psline(71.36,-34.65)(81.98,-31.11)
\psline(-9.98,-55.11)(62.02,-31.11)
\psline(0,-100)(72,-76)
\psline(0,-70)(72,-46)

\psline{->}(0,-150)(0,-170)
\psline(-120,-185)(120,-185)
\pscircle*(0,-185){5}

\end{pspicture}
}

\end{pspicture}

\caption{The blowup of \(X^{(0)} = X\) along \(E_0 = E\)}\label{F1}
\end{figure}
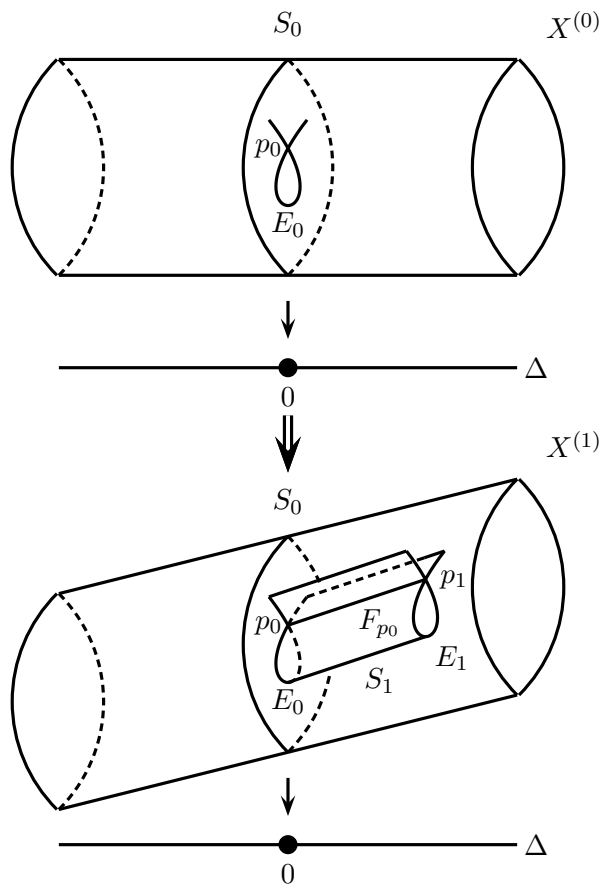

\begin{figure}[ht]
\centering
\begin{pspicture}(0,-70)(370,70)

\parametricplot[linestyle=dashed]{-6}{0}{30 0.1 t t t mul 36 sub mul mul add
5 t t mul 36 sub add}
\parametricplot[linestyle=dashed]{6}{6.84}{30 0.1 t t t mul 36 sub mul mul add
5 t t mul 36 sub add}
\parametricplot{-7.4}{-6}{30 0.1 t t t mul 36 sub mul mul add
5 t t mul 36 sub add}
\parametricplot{0}{6}{30 0.1 t t t mul 36 sub mul mul add
5 t t mul 36 sub add}
\parametricplot{6.84}{7.4}{30 0.1 t t t mul 36 sub mul mul add
5 t t mul 36 sub add}

\uput[270](30,-31){\small \(E_0\)}
\uput[180](30,5){\small \(p_0\)}
\uput[270](70,-46){\small \(S_1\)}
\uput[270](70,-10){\small \(F_{p_0}\)}
\uput[214](110,-25){\small \(p_1\)}

\psline(30,-31)(110,-61)
\psline(30,5)(110,-25)
\psline(43.882,23.76)(110,-1.034)
\psline[linestyle=dashed](110,-1.034)(123.882,-6.24)
\psline(16.118,23.76)(96.118,-6.24)

\parametricplot[linestyle=dashed]{-6}{0}{110 0.1 t t t mul 36 sub mul mul add
-25 t t mul 36 sub add}
\parametricplot[linestyle=dashed]{6}{7.4}{110 0.1 t t t mul 36 sub mul mul add
-25 t t mul 36 sub add}
\parametricplot{0}{6}{110 0.1 t t t mul 36 sub mul mul add
-25 t t mul 36 sub add}
\parametricplot{-7.4}{-6}{110 0.1 t t t mul 36 sub mul mul add
-25 t t mul 36 sub add}

\uput[270](110,-61){\small \(E_1\)}
\uput[270](150,-46){\small \(S_2\)}
\uput[270](150,-10){\small \(F_{p_1}\)}
\uput*{10pt}[230](190,5){\small \(p_2\)}

\psline(190,-31)(110,-61)
\psline(190,5)(110,-25)
\psline[linestyle=dashed](190,18.554)(123.882,-6.24)
\psline(203.882,23.76)(190,18.554)
\psline(176.118,23.76)(96.118,-6.24)

\parametricplot{6.84}{7.4}{190 0.1 t t t mul 36 sub mul mul add
5 t t mul 36 sub add}
\parametricplot[linestyle=dashed]{6}{6.84}{190 0.1 t t t mul 36 sub mul mul add
5 t t mul 36 sub add}
\parametricplot[linestyle=dashed]{-6}{0}{190 0.1 t t t mul 36 sub mul mul add
5 t t mul 36 sub add}
\parametricplot{0}{6}{190 0.1 t t t mul 36 sub mul mul add
5 t t mul 36 sub add}
\parametricplot{-7.4}{-6}{190 0.1 t t t mul 36 sub mul mul add
5 t t mul 36 sub add}

\psline(190,-31)(210,-38.5)
\psline(250,-53.5)(270,-61)

\psline(190,5)(210,-2.5)
\psline(250,-17.5)(270,-25)

\pscircle*(218,-23.5){2}
\pscircle*(230,-28.5){2}
\pscircle*(242,-33.5){2}

\psline(203.882,23.76)(223.882,16.26)

\psline(255.882,4.26)(270,-1.365)
\psline[linestyle=dashed](270,-1.365)(283.882,-6.24)

\psline(176.118,23.76)(204.118,13.26)
\psline(236.118,1.26)(256.118,-6.24)

\parametricplot[linestyle=dashed]{6}{7.4}{270 0.1 t t t mul 36 sub mul mul add
-25 t t mul 36 sub add}
\parametricplot[linestyle=dashed]{-6}{0}{270 0.1 t t t mul 36 sub mul mul add
-25 t t mul 36 sub add}
\parametricplot{0}{6}{270 0.1 t t t mul 36 sub mul mul add
-25 t t mul 36 sub add}
\parametricplot{-7.4}{-6}{270 0.1 t t t mul 36 sub mul mul add
-25 t t mul 36 sub add}

\uput[270](310,-46){\small \(S_n\)}
\uput[290](310,-10){\small \(F_{p_{n-1}}\)}
\uput[0](350,5){\small \(p_n\)}
\uput[270](270,-61){\small \(E_{n-1}\)}
\uput[270](350,-31){\small \(E_n\)}

\psline(350,-31)(270,-61)
\psline(350,5)(270,-25)
\psline[linestyle=dashed](342.623,15.786)(283.882,-6.24)
\psline(363.882,23.76)(342.623,15.786)
\psline(336.118,23.76)(256.118,-6.24)

\parametricplot{-7.4}{7.4}{350 0.1 t t t mul 36 sub mul mul add
5 t t mul 36 sub add}


\uput[90](30,60){\small \(S_0\)}
\psplot{0}{52.595}{-1 3600 4 x 30 sub x 30 sub mul mul sub sqrt mul}
\psplot{0}{58.579}{3600 4 x 30 sub x 30 sub mul mul sub sqrt}
\end{pspicture}
\caption{The blowup sequence}\label{F2}
\end{figure}
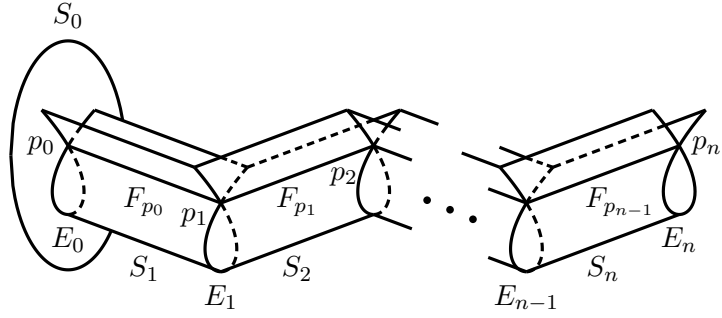

Let $E_1$ be the curve in the pencil
$|\CO_{S_1}(E_0)|$ passing through $p_1$. We blow up $X^{(1)}$ along
$E_1$ to obtain $X^{(2)}$. Now the central fiber $X_0^{(2)} = S_0 \cup
S_1 \cup S_2$ contains another ruled surface $S_2$.
Notice that we still have the exact sequence
\begin{equation}\label{s2:e7.5}
0\xrightarrow{}\CO_{E_1}\xrightarrow{} 
N_{E_1/X^{(1)}}\xrightarrow{} \CO_{E_1}\xrightarrow{} 0
\end{equation}
and hence $S_2$ is one of two ruled surfaces over $E_1\isom E$ given
in \defnref{defn1}; this is actually
true throughout our construction. If $S_2$ is twisted, we stop at
$X^{(2)}$. Otherwise, we do the same thing to $X^{(2)}$ as we did to
$X^{(1)}$. Let $F_{p_1}\subset S_2$ be the fiber of $S_2\to E_1$
over $p_1$. Notice that $X^{(2)}$ is now singular along
$F_{p_1}$, it is locally given by the equation $xy = t^2$ at a general
point of $F_{p_1}$ and there is a point $p_2\ne p_1\in F_{p_1}$ where
$X^{(2)}$ is locally given by $xy = t^2 z$. Following the convention
in \cite{C2}, we will slightly abuse the terminology to
call a singularity of the type  $xy = t^n z$ ($n > 0$) a rational
double point. Let $E_2$ be the curve in the pencil
$|\CO_{S_2}(E_1)|$ passing through the rational double point $p_2$ and
a further blowup of $X^{(2)}$ along $E_2$ will yield $X^{(3)}$.
We can continue this process and obtain a blowup sequence
\begin{equation}\label{s2:e8}
...\xrightarrow{} X^{(n)} \xrightarrow{} X^{(n-1)} \xrightarrow{}
...\xrightarrow{} X^{(1)} \xrightarrow{} X^{(0)} = X
\end{equation}
where $X_0^{(n)} = S_0\cup S_1\cup ...\cup S_n$, $S_i\cap S_{i+1} =
E_i$, $E_i\isom E$, $E_i\cdot E_{i+1} = 0$ and $S_k \isom \P^1\times
E_{k-1}$ for $1\le k\le n-1$. Let $F_{p_{n-1}}$ be the fiber of
$S_n\to E_{n-1}$ over $p_{n-1}$. \figref{F2} shows what
happens on the central fiber.

Maybe a better way to understand the singularities of the blowups
is to work out the local analytic equations of \(X^{(n)}\) over \(p\).

\begin{lem}\label{lem0.1}
Let \eqref{s2:e8} be the blowup sequence constructed as
above. Then for each \(n\ge 1\), \(X^{(n)}\) is singular along
\(F_{p_1}\cup F_{p_2}\cup ...\cup F_{p_{n-1}}\cup \{p_n\}\). At a
point \(b\in F_{p_k}\) and \(b\ne p_{k}, p_{k+1}\) for \(0\le k \le
n-1\), \(X^{(n)}\) is locally given by
\begin{equation}\label{lem0.1:e1}
xy = t^{k+1}.
\end{equation}
Locally at \(p_k\) for \(0\le k\le n-1\), 
\begin{equation}\label{lem0.1:e2}
X^{(n)}\isom \Delta_{xyzwt}^5/(xy = t^k z, zw = t)
\end{equation}
and at \(p_n\),
\begin{equation}\label{lem0.1:e3}
X^{(n)}\isom \Delta_{xyzt}^4/(xy = t^n z),
\end{equation}
where \(\Delta_{xyzwt}^5\) and \(\Delta_{xyzt}^4\) are the polydisks
parameterized by \((x, y, z, w, t)\) and \((x, y, z, t)\),
respectively.
\end{lem}

\begin{proof}
We start with \(X = X^{(0)}\) which is smooth at \(p = p_0\). Choose
local coordinates such that \(E = E_0\) is cut out by \(xy = t = 0\)
at \(p\). Blow up \(X^{(0)}\) along \(E_0\) and we obtain that
\begin{equation}\label{lem0.1:e4}
X^{(1)} \isom \Delta_{xyz_0 t}^3/(xy = tz_0),
\end{equation}
where \(z_0\) is the affine coordinate of \(F_{p_0} \isom \P^1\) such that
\(p_1\in F_{p_0}\) is given by \(z_0 = 0\) and \(p_0\) is given by
\(z_0 = \infty\). We see from
\eqref{lem0.1:e4} that \(X^{(1)}\) has a rational double point
at \(p_1\). At a point \(b\in F_{p_0}\) and \(b\ne p_0, p_1\), i.e.,
for \(z_0\ne 0,\infty\), \(X^{(1)}\) is analytically equivalent to 
\eqref{lem0.1:e1} for \(k = 0\).
At \(p_0\), i.e., at \(z_0 = \infty\), \(X^{(1)}\) is given by
\begin{equation}\label{lem0.1:e5}
xy w_0 = t
\end{equation}
where \(w_0 = 1/z_0\); this is equivalent to \eqref{lem0.1:e2} for \(k
= 0\).

Notice that \(E_1\) is cut out by \(z_0 = t = 0\). Blow up \(X^{(1)}\)
along \(E_1\) and we obtain that
\begin{equation}\label{lem0.1:e6}
X^{(2)} \isom \Delta_{xyz_1 t}^3/(xy = t^2 z_1),
\end{equation}
where \(z_1 = z_0/t\) is the affine coordinate of \(F_{p_1}\) such that
\(p_2\in F_{p_1}\) is given by \(z_1 = 0\) and \(p_1\) is given by
\(z_1=\infty\). Obviously, \(X^{(2)}\) is given by \eqref{lem0.1:e6}
at \(p_2\). At a point \(b\in F_{p_1}\) and \(b\ne p_1, p_2\), i.e.,
for \(z_1\ne 0,\infty\), \(X^{(2)}\) is analytically equivalent to 
\eqref{lem0.1:e1} for \(k = 1\).
At \(p_1\), i.e., at \(z_1 = \infty\), \(X^{(2)}\) is given by
\begin{equation}\label{lem0.1:e7}
xy = t z_0 \text{ and } w_1 z_0 = t
\end{equation}
where \(w_1 = 1/z_1\); this is equivalent to \eqref{lem0.1:e2} for
\(k = 1\).

Apply this argument inductively for \(n\) and we are done.
\end{proof}

As we will see later, the rational double point $p_n$ of
$X^{(n)}$ will play an important role in our argument.

The sequence ends at $X^{(n)}$ if $S_n\not\isom
\P^1\times E_{n-1}$ is twisted. Otherwise, let $E_n$ be the curve in
$|\CO_{S_n}(E_{n-1})|$ passing through $p_n$ and we continue to blow
up $X^{(n)}$ along $E_n$.

Suppose that $X$ is obtained from a family of K3 surfaces with a
general (and hence nonvanishing) Kodaira-Spencer class by a base
change of degree $\alpha$. We claim that the above sequence will
eventually end and it will end right at \(X^{(\alpha)}\). Namely, the
blowup sequence will end up as
\begin{equation}\label{s2:e9}
X^{(\alpha)}\xrightarrow{} X^{(\alpha-1)} \xrightarrow{}
... \xrightarrow{} X^{(1)} \xrightarrow{} X^{(0)} = X
\end{equation}
where the corresponding \(S_\alpha\subset X_0^{(\alpha)}\) is twisted.
This is clear if we reverse the process of base change and
blowups. That is, if we blow up $X$ along $E$ before we make a base
change, we will obtain $S_\alpha = \P N_{E/X}$ as the exceptional
divisor on the central fiber with indecomposable normal bundle
$N_{E/X}$ by \propref{prop1}. If we make a base change of
degree $\alpha$ afterwards, the total family $\wt{X}$ will become
singular along $E$: at a smooth point of $E$, $\wt{X}$ is locally given by
the equation $xy = t^\alpha$. We may resolve the generic singularities
of $\wt{X}$ along $E$ in the same way as in \cite[Appendix C,
p. 39]{G-H} and we will obtain a chain of ruled surfaces $S_1, S_2,
..., S_{\alpha-1}$ between $S_0 = S$ and $S_\alpha$. The resulting
family is exactly $X^{(\alpha)}$ in \eqref{s2:e9} with the required
properties.

For each \(1\le n\le \alpha\),
let $Y^{(n)}$ be the proper transform of $Y = Y^{(0)}$ under the
map $X^{(n)}\to X$. Depending on our choice of $\alpha$, the central
fiber $Y_0^{(n)}$ could be very ``bad''; for example, $Y_0^{(n)}$
could contain one or more of the double curves $E_i$ for $1\le i\le
n-1$.
However, we will show that it is possible to choose a suitable
\(\alpha\) such that the central fiber
\(Y_0^{(n)}\) of \(Y^{(n)}\) is reasonably ``well-behaved''. Most
important of all, we want to make sure that \(E_i\not\subset
Y_0^{(n)}\).

Actually, the following general statement is true, as a consequence of
the stable reduction theorem \cite{KKMS}.

\begin{thm}\label{t2}
Let \(X\) be a flat family of schemes over \(\Delta\) whose general
fibers are smooth and let \(Y\subset X\) be a closed subscheme of \(X\) of
codimension 1 which is flat over \(\Delta\). Then there exists a base
change of \(X\) followed by a series of blowups with resulting family
\(\wt{X}\) such that the proper transform \(\wt{Y}\) of \(Y\) meets
the singular locus of \(\wt{X}_0\) properly.
\end{thm}

If \(\dim X = 2\), one may think of \((X, Y)\) as a family of curves
with marked points; it is well known that after a suitable semi-stable
reduction \(\wt{X}\to X\), \(\wt{Y}\) extends to the sections of
\(\wt{X} \to \Delta\) and the marked points \(\wt{Y}_0\) can be kept
away from the singular locus of \(\wt{X}_0\). The above theorem is the
higher-dimensional analogue, which is not any harder to prove in
principle. However, we do not really need \thmref{t2} since it does
not give us any control of \(\wt{X}\) and hence cannot
be applied to our situation directly. Instead, we need the following
more precise statement.

\begin{prop}\label{prop2}
Let $X$ be a smooth family of K3 surfaces over the disk $\Delta$ whose
central fiber $X_0 = S$ is a BL K3 surface. Suppose that $X$ is
obtained from a family of K3 surfaces with a general Kodaira-Spencer
class by a base change of degree $\alpha$. Let $Y\subset X$ be a flat
family of rational curves with $Y_0\in |\CO_S(C + gF)|$ and let $E$ be
one of the 24 nodal curves in $|\CO_S(F)|$ and $m$ be the multiplicity
of $E\subset Y_0$.

Let \eqref{s2:e9} be the blowup sequence constructed as
above. Correspondingly, for each \(0\le n\le \alpha\),
let \(S_n, E_n, p_n, F_{p_n}, Y^{(n)}\) be defined as above.

Let $q_0 = C\cap E_0$ be the intersection between $C$ and $E_0$ on
$S_0$ and let $F_{q_0}\subset S_1$ be the fiber of $S_1\to E_0$ over $q_0$;
$q_i$ and $F_{q_i}$ are recursively given by letting $q_i = F_{q_{i-1}}\cap
E_i$ and $F_{q_i}\subset S_{i+1}$
be the fiber of $S_{i+1}\to E_i$ over $q_i$.

There exists a suitable choice of $\alpha$ such that the following
holds for each \(0\le n\le \alpha\):
\begin{enumerate}
\item the central fiber $Y_0^{(n)}$ of $Y^{(n)}$ does not contain
$E_i$ for $0\le i \le n - 1$;
\item $Y^{(n)}\cap S_i$ is a curve in the linear series
\begin{equation}\label{s2:e9.1}
\P H^0\left(\CO_{S_i}(m_i E_{i-1} + F_{q_{i-1}})\right)
\end{equation}
for $1\le i \le n-1$, where $m_1, m_2, ..., m_\alpha$ 
are $\alpha$ nonnegative integers satisfying $\sum_{i=1}^\alpha m_i = m$;
\item $Y^{(n)}\cap S_n = D \cup \mu E_n$, where $D$ is a curve
in the linear series
\begin{equation}\label{s2:e9.2}
\P H^0\left(\CO_{S_n}(m_n E_{n-1} + F_{q_{n-1}})\right)
\end{equation}
and $\mu = \sum_{i=n+1}^\alpha m_i$;
\item $F_{q_i}\subset \left(Y^{(n)}\cap S_i\right)$ for
$1\le i\le n \le \alpha - 1$.
\end{enumerate}
\end{prop}

Although the general results on stable reduction such as \thmref{t2}
cannot be applied to \propref{prop2} directly, its proof is actually
carried out by explicitly applying semi-stable reduction to
\(X^{(\alpha)}\).

By \propref{prop2}, $Y_0^{(\alpha)}$ looks as follows:
the components of $Y_0^{(\alpha)}$ over $E$ consist of
\begin{equation}\label{s2:e9.3}
(F_{q_0}\cup D_1) \cup (F_{q_1}\cup D_2)\cup
... \cup (F_{q_{\alpha - 2}} \cup D_{\alpha - 1}) \cup \Gamma,
\end{equation}
where
$D_i \subset S_i$, $D_i\in |\CO_{S_i}(m_i E_{i-1})|$, $E_i\not\subset
D_i$ for $1\le i \le \alpha - 1$, $\Gamma\subset S_\alpha$ and $\Gamma
\in |\CO_{S_\alpha}(m_{\alpha} E_{\alpha - 1} + F_{q_{\alpha-1}})|$.
We will call the components $D_i$ ``wandering components''. Actually,
we have

\begin{prop}\label{prop3}
With all the notations as above, then
\begin{equation}\label{prop3:e1}
m_1 = m_2 = ... = m_{\alpha - 1} = 0,
\end{equation}
i.e., \(D_i = \emptyset\) for \(1\le i\le \alpha - 1\) and
there are no wandering components at all.
\end{prop}

Therefore, ``interesting'' things only happen on the twisted ruled
surface $S_\alpha$. Among the components $F_{q_0}\cup F_{q_1}\cup
... \cup F_{q_{\alpha - 2}} \cup \Gamma$ of $Y_0^{(\alpha)}$,
$F_{q_i}$'s are a chain of rational curves connecting $C$ and
$\Gamma$ and they will be contracted under stable reduction; the only
nontrivial part is $\Gamma\subset S_\alpha$ which maps to $E$ with a
degree $m$ map. One of main steps of our proof is to classify all possible
configurations of $\Gamma$.

Let $\delta(A)$ denote the total $\delta$-invariant of the
singularities of a curve $A$ and let $\delta(A, B)$ denote the total
$\delta$-invariant of the singularities of $A$ in the (analytic) neighborhood of
$B$. The latter notation $\delta(A, B)$ is used under two circumstances:
\begin{enumerate}
\item if $B\subset A$ is a closed subscheme of $A$,
$\delta(A, B)$ is simply the total
$\delta$-invariant of the singularities $p$ of $A$ with $p\in B$;
\item if $\Upsilon$ is a family of curves over the disk $\Delta$,
$A = \Upsilon_t$ is the general fiber of $\Upsilon\to\Delta$ and
$B\subset \Upsilon_0$ is a closed subscheme of the central fiber
$\Upsilon_0$,
then $\delta(\Upsilon_t, B)$ is the total $\delta$-invariant of the
singularities of $\Upsilon_t$ in the neighborhood of $B$; notice that this
is well-defined.
\end{enumerate}

We claim that

\begin{prop}\label{prop4}
Suppose that \propref{prop2} and \ref{prop3} are true.
With all the notations as above,
\begin{equation}\label{s2:e10}
\delta(Y_t^{(\alpha)}, \Gamma)\ge m
\end{equation}
and if the equality holds, the general fiber $Y_t^{(\alpha)}$ of
$Y^{(\alpha)}$ has exactly $m$ nodes in the neighborhood of $\Gamma$.
Or equivalently, $\delta(Y_t, E) \ge m$
and if the equality holds, the general fiber $Y_t$ of $Y$ has exactly $m$
nodes in the neighborhood of $E$.
\end{prop}

Notice that the total $\delta$-invariant of $Y_t$ is $g$ and
\begin{equation}\label{s2:e11}
g = \delta(Y_t) \ge \sum \delta(Y_t, E)
\end{equation}
where we sum over all the $24$ nodal fibers $E$ of $S\to\P^1$. By
\propref{prop4}, the RHS of \eqref{s2:e11} is at least the sum of the
multiplicities of $E$ in $Y_0$, which is $g$. So we must have
$\delta(Y_t, E) = m_E$ for each nodal fiber $E$, where $m_E$ is the
multiplicity of $E$ in $Y_0$. By \propref{prop4} again, $Y_t$ is nodal
in the neighborhood of each $E$. And our main theorem follows.

The rest of the paper is organized as follows. In \secref{s3}, we will
introduce some preliminary results that will be needed later in our
proof, which include a geometrical
construction of the twisted ruled surface over $E$
and some local results on the deformation of curve singularities.
Next we will prove \propref{prop4} in \secref{s4}, during which we
will give a classification for all possible configurations of $\Gamma$ and
the stable reduction over it. The proofs of \propref{prop2} and
\ref{prop3} will be postponed until \secref{s5}.

\section{Preliminaries}\label{s3}

\subsection{Construction of the Twisted Ruled Surface}\label{s3:1}

Let $E$ be a rational curve with one node and let $W$ be a rank 2
vector bundle over $E$ satisfying the exact sequence \eqref{defn1:e1}.
As mentioned before, there are two isomorphism classes of $\P W$:
one is ``trivial'' and the other is ``twisted''.
We will give an explicit geometric construction of the latter.

Let $\nu: \wt{E}\to E$ be the normalization of $E$. Since $\wt{E}\isom
\P^1$, $\nu^*(W)$ splits to $\CO_{\wt{E}} \oplus \CO_{\wt{E}}$ on
$\wt{E}$. And this induces a map $\nu: \P^1\times \wt{E} \isom
\P^1\times \P^1 \to \P W$, which is just the normalization of $\P
W$. Intuitively, we say $\nu$ ``unfolds'' $\P W$.

We use $E$ and $\wt{E}$ to denote the zero sections of
$\P W\to E$ and its normalization \(\P^1\times \wt{E}\to \wt{E}\),
respectively.

Let $a, b\in \wt{E}$ be the preimages of the node $p\in E$. Let $F_a,
F_b$ be the fibers of $\P^1\times \wt{E}$ over $a, b$ and let $F_p$ be
the fiber of $\P W\to E$ over $p$. One can think of $\P W$ being constructed
from $\P^1\times \wt{E}$ by ``gluing'' two fibers $F_a$ and $F_b$.

Let $\nu_a: F_a \to F_p$ and $\nu_b: F_b \to F_p$ be the maps induced
by $\nu$. We have a natural identification $\phi_{ab}$ between $F_a$
and $F_b$ on $\P^1\times \wt{E}$, which simply sends $x\in F_a$ to
$y\in F_b$ if there is a curve in the pencil $|\wt{E}|$ passing
through $x$ and $y$. So $h = \phi_{ba} \circ \nu_b^{-1} \circ \nu_a$
is an automorphism of $F_a\isom \P^1$, where $\phi_{ba} =
\phi_{ab}^{-1}$.

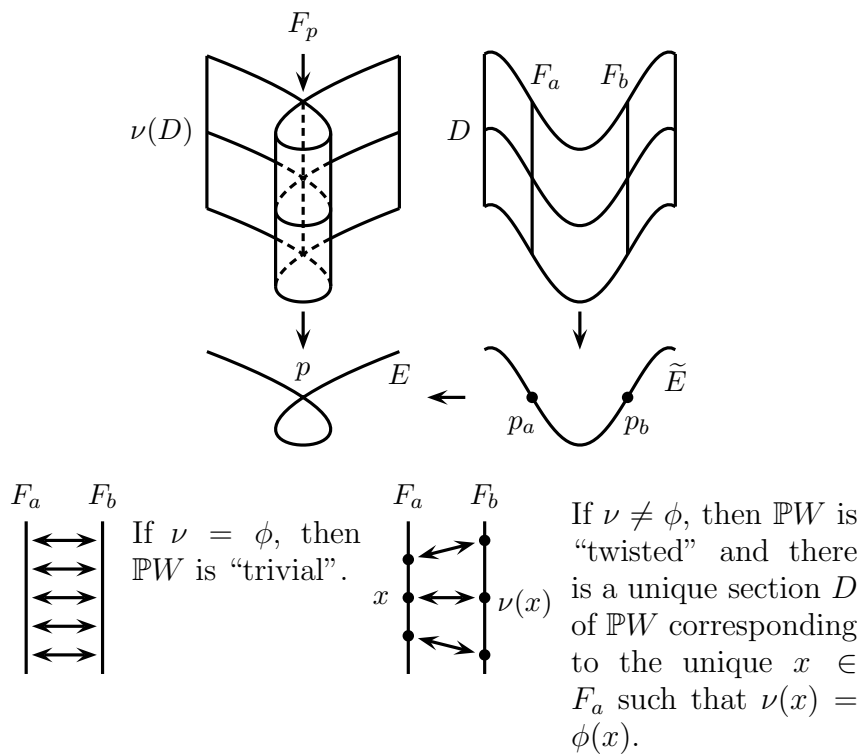
\begin{figure}[ht]
\centering
\begin{pspicture}(-10,-370)(430,20)

\rput[lt]{0}(70,20){
\begin{pspicture}(-10,-210)(105,20)
\uput[90](55,-5){\small \(F_p\)}

\psline{->}(55,-5)(55,-25)

\psline[linestyle=dashed](55,-30)(55,-110)

\parametricplot{-7}{7}{55 0.3 t t t mul 25 sub mul mul add
t t mul 25 sub 30 sub}

\parametricplot{5.774}{7}{55 0.3 t t t mul 25 sub mul mul add
t t mul 25 sub 70 sub}
\parametricplot{-7}{-5.774}{55 0.3 t t t mul 25 sub mul mul add
t t mul 25 sub 70 sub}
\parametricplot{-2.887}{2.887}{55 0.3 t t t mul 25 sub mul mul add
t t mul 25 sub 70 sub}
\parametricplot[linestyle=dashed]{2.887}{5.774}{%
55 0.3 t t t mul 25 sub mul mul add
t t mul 25 sub 70 sub}
\parametricplot[linestyle=dashed]{-5.774}{-2.887}{%
55 0.3 t t t mul 25 sub mul mul add
t t mul 25 sub 70 sub}

\uput[180](4.6,-46){\small \(\nu(D)\)}

\parametricplot{5.774}{7}{55 0.3 t t t mul 25 sub mul mul add
t t mul 25 sub 110 sub}
\parametricplot{-7}{-5.774}{55 0.3 t t t mul 25 sub mul mul add
t t mul 25 sub 110 sub}
\parametricplot{-2.887}{2.887}{55 0.3 t t t mul 25 sub mul mul add
t t mul 25 sub 110 sub}
\parametricplot[linestyle=dashed]{2.887}{5.774}{%
55 0.3 t t t mul 25 sub mul mul add
t t mul 25 sub 110 sub}
\parametricplot[linestyle=dashed]{-5.774}{-2.887}{%
55 0.3 t t t mul 25 sub mul mul add
t t mul 25 sub 110 sub}

\psline(40.566,-46.667)(40.566,-126.667)
\psline(69.434,-46.667)(69.434,-126.667)

\psline(4.6,-6)(4.6,-86)
\psline(105.4,-6)(105.4,-86)

\psline{->}(55,-140)(55,-160)
\uput[270](55,-160){\small \(p\)}

\parametricplot{-7}{7}{55 0.3 t t t mul 25 sub mul mul add
t t mul 25 sub 185 sub}

\uput[270](105.4,-161){\small \(E\)}
\end{pspicture}
}

\rput[lt]{0}(230,0){
\begin{pspicture}(0,-210)(100,0)
\parametricplot{-10}{10}{50 5 t mul add
-30 t t mul 25 sub 150 t t mul sub mul 150 div add}
\parametricplot{-10}{10}{50 5 t mul add
-70 t t mul 25 sub 150 t t mul sub mul 150 div add}
\parametricplot{-10}{10}{50 5 t mul add
-110 t t mul 25 sub 150 t t mul sub mul 150 div add}

\psline(0,-5)(0,-85)
\psline(100,-5)(100,-85)

\uput[65](25,-30){\small \(F_a\)}
\uput[115](75,-30){\small \(F_b\)}

\uput[180](0,-45){\small \(D\)}

\psline(25,-30)(25,-110)
\psline(75,-30)(75,-110)

\psline{->}(50,-140)(50,-160)

\parametricplot{-10}{10}{50 5 t mul add
-185 t t mul 25 sub 150 t t mul sub mul 150 div add}

\psline{<-}(-30,-185)(-10,-185)

\pscircle*(25,-185){3}
\pscircle*(75,-185){3}

\uput[250](25,-188){\small \(p_a\)}
\uput[290](75,-188){\small \(p_b\)}

\uput[270](100,-160){\small \(\wt{E}\)}
\end{pspicture}
}

\rput[lt]{0}(-10,-230){
\begin{pspicture}(0,-100)(170, 0)

\uput[90](0,-20){\small \(F_a\)}
\uput[90](40,-20){\small \(F_b\)}

\psline(0,-20)(0,-100)
\psline(40,-20)(40,-100)
\psline{<->}(3,-30)(37,-30)
\psline{<->}(3,-45)(37,-45)
\psline{<->}(3,-60)(37,-60)
\psline{<->}(3,-75)(37,-75)
\psline{<->}(3,-90)(37,-90)

\rput[lt]{0}(50,-20){
\parbox{1.2in}{If \(\nu = \phi\), then \(\P W\) is ``trivial''.}
}
\end{pspicture}
}

\rput[lt]{0}(190,-230){
\begin{pspicture}(0,-100)(230,0)

\uput[90](0,-20){\small \(F_a\)}
\uput[90](40,-20){\small \(F_b\)}

\psline(0,-20)(0,-100)
\psline(40,-20)(40,-100)

\psline{<->}(5,-60)(35,-60)
\psline{<->}(5,-38.75)(35,-31.25)
\psline{<->}(5,-81.25)(35,-88.75)

\pscircle*(0,-40){3}
\pscircle*(40,-30){3}
\pscircle*(0,-60){3}
\pscircle*(40,-60){3}
\pscircle*(0,-80){3}
\pscircle*(40,-90){3}

\uput[180](-3,-60){\small \(x\)}
\uput[0](40,-63){\small \(\nu(x)\)}

\rput[lt]{0}(80,-10){
\parbox{1.5in}{If \(\nu \ne \phi\), then \(\P W\) is ``twisted''
and there is a unique section \(D\) of \(\P W\) corresponding to the
unique \(x\in F_a\) such that \(\nu(x) = \phi(x)\).}}
\end{pspicture}
}
\end{pspicture}

\caption{``trivial'' vs ``twisted''}\label{F3}
\end{figure}

If $x\in F_a$ is a fixed point of $h$, i.e., $h(x) = x$, the curve
$D\in |\wt{E}|$ passing through $x$ and $\phi_{ab}(x)$ maps to a curve
$\nu(D)\in |E|$.
If $W$ is indecomposable, there is only one curve in $|E|$. So
$h$ can have only one fixed point. If we represent $h$ by a matrix \(H\in
GL(2)\), $H$ has only one eigenvector and is hence equivalent to
\begin{equation}\label{s3:e1}
\begin{pmatrix}
1 & \lambda\\
0 & 1
\end{pmatrix}.
\end{equation}
In fact, $\lambda$ in \eqref{s3:e1} classifies all the extensions in
$\Ext(\CO_E, \CO_E) = \BC$. For $\lambda = 0$, we obtain $\P^1\times E$; for
$\lambda\ne 0$, we obtain $\P W$ with $W$ indecomposable and they are
isomorphic to each other (see \figref{F3}).

\begin{rem}\label{s3:rem1}
On an interesting though unrelated issue, one may ask what kind of
surfaces we get if we glue \(\P^1\times \P^1\) along \(F_a\) and
\(F_b\) via an automorphism \(h\) which has two fixed points, i.e.,
whose corresponding matrix representation \(H\) has two
eigenvectors. The resulting surface \(S\) will have exactly two sections
\(D_1\) and \(D_2\) with self-intersection zero. So what kind of
surface is \(S\)? Actually, \(S = \P (\CO_E \oplus L_E)\) where \(L_E\) is
a nontrivial line bundle on \(E\) with \(\deg L_E = 0\). The two
sections \(D_1\) and \(D_2\) are not linearly equivalent on \(S\) and
they correspond to the global sections of \(\CO_E \oplus L_E\) and
\(\CO_E \oplus L_E^{-1}\), respectively. I would like to thank James
McKernan for pointing this out to me.
\end{rem}

\subsection{A Key Lemma}

This is basically Lemma 2.2 in \cite{C1} or Lemma 2.1 in \cite{C2}.

\begin{lem}\label{lem1}
Let $X\subset \Delta_{xyz}^3\times \Delta_t$ be a family of surfaces
given by $xy = t^\alpha$ for some $\alpha > 0$. 
Let $X_0$ be the central fiber of $X$ over $\Delta_t$ and
$X_0 = R_1 \cup R_2$ where $R_1 = \{x = t = 0\}$
and $R_2 = \{y = t = 0\}$ and let $E = R_1 \cap R_2$.
Let $Y$ be a flat family of curves over
$\Delta_t$ and $\pi: Y\to X$ be a proper morphism preserving
the base $\Delta_t$. Suppose that $E\not\subset \pi(Y_0)$, where $Y_0$ is
the central fiber of $Y$.
Let $Y_0 = \Gamma_1 \cup\Gamma_2$ with $\pi(\Gamma_1)\subset R_1$ and
$\pi(\Gamma_2)\subset R_2$. Then $\pi(\Gamma_1) \cdot E =
\pi(\Gamma_2) \cdot E$, where the intersections $\pi(\Gamma_1)\cdot E$
and $\pi(\Gamma_2)\cdot E$ are taken on the surfaces $R_1$ and $R_2$,
respectively.
\end{lem}

The proof of this lemma is not hard. The readers may find a proof in
\cite{C1} or \cite{C2}.


\begin{figure}[ht]
\centering
\begin{pspicture}(-10,-270)(380,30)

\psline[linestyle=dashed](0,-70)(300,-70)
\psline(300,-70)(300,-190)(0,-190)(0,-70)
\psline(0,-70)(70,0)
\psline[linestyle=dashed](70,0)(370,0)
\psline(370,0)(300,-70)
\psline(370,0)(370,-120)
\psline(370,-120)(300,-190)

\psline{->}(150,-210)(150,-230)

\pscircle*(150,-250){5}

\psline(0,-250)(300,-250)

\uput[270](150,-260){\small \(0\)}
\uput[0](300,-250){\small \(\Delta\)}

\uput[45](370,0){\small \(X\)}

\psline(35,-35)(120,-35)
\psline(250,-35)(335,-35)
\psline[linestyle=dashed](120,-35)(250,-35)

\pscurve(120,-35)(105,-105)(120,-155)
\pscurve(90,-35)(80,-105)(90,-155)
\pscurve(60,-35)(55,-105)(60,-155)

\psline(35,-35)(35,-155)

\psline(335,-35)(335,-155)

\psline(35,-155)(120,-155)
\psline(250,-155)(335,-155)
\psline[linestyle=dashed](120,-155)(250,-155)

\pscurve(250,-35)(235,-105)(250,-155)
\pscurve(280,-35)(270,-105)(280,-155)

\psline(0,-190)(35,-155)

\psline(150,-190)(220,-120)
\psline(150,-70)(220,0)

\psline(130,-140)(200,-70)
\uput[315](140,-130){\small \(E\)}

\pscurve(130,-140)(140,-155)(150,-190)
\pscurve(150,-70)(140,-115)(130,-140)

\uput[135](185,-35){\small \(R_1\)}
\uput[315](185,-155){\small \(R_2\)}

\pscurve(165,-105)(175,-120)(185,-155)
\pscurve(185,-35)(175,-80)(165,-105)

\uput[21](175,-130){\small \(\Gamma_2\)}
\uput[344](180,-50){\small \(\Gamma_1\)}

\pscurve(200,-70)(210,-85)(220,-120)
\pscurve(220,0)(210,-45)(200,-70)

\end{pspicture}
\caption{\lemref{lem1}}\label{F4}
\end{figure}

\begin{defn}\label{defn2}
Let $Y$ be a one-parameter family of curves over $\Delta$ and let $p\in
Y_0$ be a point on the central fiber $Y_0$. Even if $Y$ is irreducible
globally, it is still possible that $Y$ is reducible in an analytic
neighborhood of $p$.
That is, if we let $U$ be an analytic neighborhood of $Y$ at $p$, $U$
might be reducible such that $U = \cup V_i$ where we call each $V_i$
a {\it local irreducible component\/} of $Y$ at $p$.
This happens if $Y$ is not normal and the general
fiber $Y_t$ is singular in the neighborhood of $p$.
If $Y$ breaks
into several local irreducible components at $p$, the normalization of
$Y$ will make these components disconnected. Let $\Gamma_1$ and
$\Gamma_2$ be two local branches of $Y_0$ at $p$. We call $\Gamma_1$
is {\it locally separated from\/} $\Gamma_2$ at $p$ if 
$\Gamma_1$ and $\Gamma_2$ do not lie on the same local irreducible
component of $Y$ at $p$, or
equivalently, $\Gamma_1$ and $\Gamma_2$ become disconnected on the
normalization of $Y$. And we call $Y$ is {\it totally separated\/} at
$p$ if any two branches of $Y_0$ at $p$ are locally separated from each
other, i.e., if $Y_0 = \mu_1 \Gamma_1 \cup \mu_2 \Gamma_2 \cup
... \cup \mu_k \Gamma_k$ at $p$, $\Gamma_i$ is locally
separated from $\Gamma_j$ for all $1\le i\ne j\le k$.
\end{defn}

\begin{rem}\label{rem2}
It is necessary to point out that \lemref{lem1} is a local result. So
it holds for every local irreducible component of $Y$ at
$\pi^{-1}(p)$, where $p\in X$ is the origin.
For example, suppose that $Y_0 = \Gamma_1\cup \Gamma_2$ with
$\pi(\Gamma_i)\subset R_i$ for $i = 1,2$ and $\Gamma_1$ is
reduced and locally
irreducible. Then we certainly have $\pi(\Gamma_1) \cdot E =
\pi(\Gamma_2)\cdot E$ by the lemma; in particular, this means
$\Gamma_2\ne \emptyset$. In addition, we can also conclude by the lemma
that $Y$ is locally irreducible at $\pi^{-1}(p)$, which implies that no
component of $\Gamma_2$  is locally separated from $\Gamma_1$.
As for another example, take $Y_0 = \cup_{i=1}^4 \Gamma_i$ with
$\pi(\Gamma_1), \pi(\Gamma_2)\subset R_1$ and
$\pi(\Gamma_3), \pi(\Gamma_4)\subset R_2$ and suppose that each
$\pi(\Gamma_i)$ meets $E$ transversely.
Then we may conclude by the lemma that $Y$ consists of at most two
local irreducible components and if this happens, we have either
$\Gamma_1$ and $\Gamma_3$ lie on one component and $\Gamma_2$ and
$\Gamma_4$ lie on the other or $\Gamma_1$ and $\Gamma_4$
lie on one component and $\Gamma_2$ and
$\Gamma_3$ lie on the other; in particular, $Y$ cannot be totally
separated at $\pi^{-1}(p)$.
\end{rem}

For a three-fold rational double point $p\in X$ given by
$xy = t^\alpha z$, we
can resolve $X$ at $p$ by blowing up one of the two surfaces of
$X_0$ at $p$, i.e., let $\wt{X}\subset X\times \P^1$ be the resolution
given by
\begin{equation}\label{s3:e2}
\frac{x}{z} = \frac{t^\alpha}{y} = \frac{W_1}{W_0},
\end{equation}
where $(W_0, W_1)$ is the homogeneous coordinate of $\P^1$.
Strictly speaking, it is not a resolution of singularities because
$\wt{X}$ is still singular if $\alpha > 1$. But now $\wt{X}$ is given
by $wy = t^\alpha$ along its singular locus, where we may apply
\lemref{lem1} to obtain the following corollary.

\begin{cor}\label{cor1}
Let $X, R_1, R_2, E, \pi, Y$ be defined as in \lemref{lem1} except that $X$
is given by $xy = t^\alpha z$ instead. Suppose that $Y_0$ contains a
component $\Gamma_1$ such that $\pi(\Gamma_1)\subset R_1$
is tangent to $E$ at the origin $p$. Then there must exist a component
$\Gamma_2$ of $Y_0$ such that $\pi(\Gamma_2)\subset R_2$ passes
through $p$. In particular, $Y$ cannot be totally separated at
point $q$ where $\pi(q) = p$ and $q\in \Gamma_1$.
\end{cor}

\begin{proof}
See \cite[Corollary 2.1]{C2}.
\end{proof}

\subsection{Some Results on Curve Singularities}

The following lemma is
basically a combination of Corollary 4.1 and Proposition 4.3
in \cite{C2}.

\begin{lem}\label{lem2}
Let $Y\subset \Delta^2\times\Delta_t$ be a reduced flat family of
curves over $\Delta_t$ with central fiber $Y_0 = \mu_1 \Gamma_1 \cup
\mu_2 \Gamma_2\cup ...\cup \mu_n \Gamma_n$, where $\mu_i$ is the
multiplicity of the component $\Gamma_i$ in $Y_0$. Suppose that $Y$ is
totally separated at the origin $p$. Then
\begin{equation}\label{s3:e3}
\delta(Y_t) \ge \sum_{1\le r < s\le n} \mu_r\mu_s(\Gamma_r\cdot
\Gamma_s)
\end{equation}
where the intersections $\Gamma_r\cdot \Gamma_s$ are taken on $\{t =
0\}\isom \Delta^2$.

If the equality holds in \eqref{s3:e3} and we further assume that
\begin{enumerate}
\item[A1.] $\Gamma_r$ and $\Gamma_s$ meet transversely, i.e.,
$\Gamma_r\cdot \Gamma_s = 1$ for $1\le r < s\le n$, and
\item[A2.] for each irreducible component $Z\subset Y$ of $Y$,
the central fiber $Z_0$ of $Z$ is reduced, i.e., $Y$ consists of
exactly $\sum_{i=1}^n \mu_i$ irreducible components,
\end{enumerate}
then $Y_t$ is nodal.
\end{lem}

\begin{proof}
See \cite[Sec. 4]{C2}.
\end{proof}

\begin{rem}\label{rem3}
Here is an example how to apply \lemref{lem2}. Let 
\begin{equation}\label{rem3:e1}
Y\subset \Delta_{xy}^2\times \Delta_t
\end{equation}
be a reduced flat family of curves whose
central fiber $Y_0$ is given by $x^m y^n = 0$, i.e., $Y_0 = m\Gamma_1
\cup n\Gamma_2$ where $\Gamma_1$ and $\Gamma_2$ are the curves $\{x =
t = 0\}$ and $\{y = t = 0\}$, respectively. Suppose that $Y$ is
totally separated at the origin $p$. That is to say that for each
irreducible component $Z\subset Y$, either $Z_0 = m'\Gamma_1$ for some
$m' \le m$ or $Z_0 = n'\Gamma_2$ for some $n' \le n$. Then
\lemref{lem2} yields that $\delta(Y_t) \ge mn$. If we further assume
that $\delta(Y_t) = mn$ and $Y$ has exactly $m+n$ irreducible
components, then $Y_t$ has exactly $mn$ nodes as singularities.
\end{rem}

The above lemma can be applied to a family of curves in the
neighborhood of a three-fold rational double point $xy = t^\alpha z$.

\begin{cor}\label{cor2}
Let $X\subset \Delta_{xyz}^3\times \Delta_t$ be a family of surfaces
given by $xy = t^\alpha z$ for some $\alpha > 0$ and let $R_1, R_2, E$ be
defined as in \lemref{lem1}. Let $Y\subset X$ be a reduced closed
subscheme of $X$ with codimension $1$ and suppose that $E\not\subset
Y_0$. Let $Y_0 = \Gamma_1\cup \Gamma_2$, where $\Gamma_1\subset R_1$
and $\Gamma_2\subset R_2$. If
\begin{enumerate}
\item[A1.] each irreducible component of $Y_0$ meets $E$ transversely and
\item[A2.] $Y$ is totally separated at the origin $p$,
\end{enumerate}
then
\begin{equation}\label{s3:e4}
\delta(Y_t) \ge \mu_1\mu_2
\end{equation}
where $\mu_1 = \Gamma_1\cdot E$ and $\mu_2 = \Gamma_2\cdot E$.
If the equality holds in \eqref{s3:e4} and we further assume that
\begin{enumerate}
\item[A3.] for each irreducible component $Z\subset Y$ of $Y$,
the central fiber $Z_0$ of $Z$ is reduced, i.e.,
$Y$ consists of exactly $\mu_1 + \mu_2$ irreducible components,
\end{enumerate}
then $Y_t$ is nodal.
\end{cor}

This is a weak version of Proposition 4.4 and 4.5 in
\cite{C2}, which can be proved by first resolving $X$ as in
\eqref{s3:e2} and then applying \lemref{lem2}. Please see \cite[Sec. 4]{C2}
for the details.

\section{Proof of \propref{prop4}}\label{s4}

First we ``unfold'' the twisted ruled surface $S_\alpha$ as in
\ssecref{s3:1}. Let 
\begin{equation}\label{s4:e-2}
\nu: \wt{S}_{\alpha} = \P^1\times\wt{E}_{\alpha-1}\to S_\alpha
\end{equation}
be the normalization of $S_\alpha$,
where $\wt{E}_{\alpha-1}$ is the normalization of $E_{\alpha-1}$. Let
$a, b\in \wt{E}_{\alpha-1}$ be the preimages of the node
$p_{\alpha-1}$ and let $F_a, F_b\subset \wt{S}_{\alpha}$ be the fibers
over $a$ and $b$. Let $\nu_a: F_a\to F_{p_{\alpha-1}}$ and $\nu_b:
F_b\to F_{p_{\alpha-1}}$ be the maps induced by $\nu$ and let
$\varepsilon_{ab} = \nu_b^{-1} \circ \nu_a$ and $\varepsilon_{ba} =
\nu_a^{-1} \circ \nu_b$. We will abbreviate both $\varepsilon_{ab}$
and $\varepsilon_{ba}$ to $\varepsilon$ most of time since it is
usually clear which one we are using, i.e., $\varepsilon(u) =
\varepsilon_{ab}(u)$ if $u\in F_a$ and $\varepsilon(u) =
\varepsilon_{ba}(u)$ if $u\in F_b$. Also we write
$u\xrightarrow{\varepsilon} w$ if $w = \varepsilon(u)$.

Let $\phi_{ab}$ and $\phi_{ba}$ be defined as in \ssecref{s3:1}, i.e.,
$w = \phi_{ab}(u)$ if $u\in F_a$ and $w\in F_b$ lie on a curve in the pencil
$|\CO_{\wt{S}_\alpha}(\wt{E}_{\alpha-1})|$. Again, we will abbreviate
both $\phi_{ab}$ and $\phi_{ba}$ to $\phi$, i.e., $\phi(u) =
\phi_{ab}(u)$ if $u\in F_a$ and $\phi(u) = \phi_{ba}(u)$ if $u\in F_b$.
We write $u\xrightarrow{\phi} w$ if $w = \phi(u)$. Also we use the
notation $\overline{uw}$ to denote the curve in
$|\CO_{\wt{S}_\alpha}(\wt{E}_{\alpha-1})|$ passing through $u$ and $w$
if $u\xrightarrow{\phi} w$.

Let $r_a\in F_a$ and $r_b\in F_b$ be the preimages of the rational
double point $p_\alpha$. Using the notations just defined, we have
$r_a\xrightarrow{\varepsilon} r_b$ and $r_b\xrightarrow{\varepsilon} r_a$.

Let $\wt{\Gamma} = \nu^{-1}(\Gamma)\subset\wt{S}_\alpha$.
Suppose that $\wt{\Gamma}$ meets
$F_a$ at a point $u\ne r_a$ with multiplicity $k$.
The branches of
$\wt{\Gamma}$ at $u$ map to the branches of $\Gamma$ lying on one of
the two surfaces of $X_0^{(\alpha)}$ at $\nu(u)$,
where $X^{(\alpha)}$ is locally given by
$xy = t^\alpha$. So we can apply \lemref{lem1} to $Y^{(\alpha)}\subset
X^{(\alpha)}$ at $\nu(u)$ and conclude that there must be branches of
$\Gamma$ lying on the other surface of $X_0^{(\alpha)}$ at $\nu(u)$
and the branches on both surfaces must meet $F_{p_{\alpha-1}}$ at
$\nu(u)$ with the same multiplicity $k$. Correspondingly, $\wt{\Gamma}$
must meet $F_b$ at $w = \varepsilon(u)$ with multiplicity
$k$. Therefore, if $\wt{\Gamma}$ meets $F_a$ at $u\ne r_a$ with
multiplicity $k$, $\wt{\Gamma}$ must meet $F_b$ at $w =
\varepsilon(u)$ with the same multiplicity $k$. Similarly, if
$\wt{\Gamma}$ meets $F_b$ at $w\ne r_b$ with multiplicity $k$,
$\wt{\Gamma}$ must meet $F_a$ at $u = \varepsilon(w)$ with the same
multiplicity $k$. So we can pair each $u\ne r_a\in \wt{\Gamma}\cap
F_a$ with $w = \varepsilon(u)\ne r_b\in \wt{\Gamma}\cap F_b$ and
$(\wt{\Gamma}\cdot F_a)_u = (\wt{\Gamma}\cdot F_b)_w$. And for the
remaining pair $r_a\xrightarrow{\varepsilon} r_b$, we must have
$(\wt{\Gamma}\cdot F_a)_{r_a} = (\wt{\Gamma}\cdot F_b)_{r_b}$.
In summary, we have
\begin{equation}\label{s4:e-1}
(\wt{\Gamma}\cdot F_a)_{u} = (\wt{\Gamma}\cdot F_b)_{w}
\end{equation}
for any pair of points $u\in F_a$ and $w\in F_b$ with
$u\xrightarrow{\varepsilon} w$.

Let $N\subset \Gamma$ be the irreducible component of
$\Gamma$ with 
\begin{equation}\label{s4:e0}
N\in |\CO_{S_\alpha}(F_{q_{\alpha-1}} + \mu E_{\alpha-1})|
\end{equation}
for some $\mu \le m$. And let $\wt{N} = \nu^{-1}(N)\subset
\wt{S}_\alpha$.

Let $\wt{Y}\to Y^{(\alpha)}$ be the stable reduction of $Y^{(\alpha)}$
after normalization. Namely, $\wt{Y}_t$ is the normalization of
$Y_t^{(\alpha)}$ on the general fibers and
\begin{equation}\label{s4:e0.1}
\wt{Y}_0\to Y_0^{(\alpha)}
\end{equation}
is a stable map on the central fiber.
We say a component $M_1\subset\wt{Y}_0$ is joined to another
component $M_2\subset \wt{Y}_0$ over a point $s\in Y_0^{(\alpha)}$ if
the two components $M_1$ and $M_2$ are joined by a chain of curves
contracted to $s$.


Consider the component of $\wt{Y}_0$ that dominates $N$. It must be
isomorphic to $\wt{N}\subset \wt{S}_\alpha$. So we use the same
notation $\wt{N}$ to denote this component.

We call a sequence of points $\{u_0, w_0, u_1, w_1, ..., u_n, w_n\}\subset
\wt{\Gamma} \cap (F_a\cup F_b)$ an {\it S-chain\/} if $u_0\in F_a$
and
\begin{equation}\label{s4:e1}
u_0\xrightarrow{\varepsilon} w_0\xrightarrow{\phi} u_1
\xrightarrow{\varepsilon} w_1\xrightarrow{\phi} ...\xrightarrow{\varepsilon}
w_{n-1}\xrightarrow{\phi} u_n\xrightarrow{\varepsilon} w_n.
\end{equation}
Notice that $u_{i+1} = h(u_i)$ where $h = \phi\circ
\varepsilon\in\Aut(F_a) \isom \Aut(\P^1)$ is the automorphism of $\P^1$
given by \eqref{s3:e1} with $\lambda\ne 0$
if we let $a\in F_a$ be the point at $\infty$.
Obviously, $h^k(u) \ne u$ for any $u\ne a$ and $k \ne 0$ and hence
$u_i \ne u_j$ for any $i\ne j$. Similarly, $w_i\ne w_j$ for any $i\ne j$.
Therefore, the points in an S-chain are distinct.

An S-chain is maximal if it is not contained in a longer S-chain. We
claim that

\begin{prop}\label{prop5}
A maximal S-chain must contain either $r_a$ or $r_b$.
\end{prop}

\begin{proof}
Let $\{u_0, w_0, u_1, w_1, ..., u_n, w_n\}$ be a maximal S-chain and
\begin{equation}\label{prop5:e1}
r_a, r_b\not\in\{u_0, w_0, u_1, w_1, ..., u_n, w_n\}.
\end{equation}



Since $\{u_0, w_0, u_1, w_1, ..., u_n, w_n\}$ is maximal, there does
not exists $w\in \wt{\Gamma}\cap F_b$
such that $w\xrightarrow{\phi} u_0$ and there is no
curve $\overline{wu_0}\subset \wt{\Gamma}$. So $\wt{N}$ has to pass
through $u_0$. Similarly, there is no point $u\in F_a$ such that
$\overline{w_n u}\subset \wt{\Gamma}$ and hence $\wt{N}$ must pass
through $w_n$.

Applying \lemref{lem1} to the point $\nu(u_0) = \nu(w_0)$, we see that
the branch of $\wt{N}$ at $u_0$ is joined to either the branch of
$\wt{N}$ at $w_0$ or a component $M_1$ dominating
$\nu(\overline{w_0u_1})$ over $\nu(u_0)$.
If it is the former case that the branch of $\wt{N}$ at $u_0$
is joined to the branch of $\wt{N}$ at $w_0$ over $\nu(u_0)$,
it contradicts the fact that the dual graph of $\wt{Y}_0$ is a tree.
Otherwise, if $\wt{N}$ is joined to $M_1$ over $\nu(u_0)$,
we continue to apply
\lemref{lem1} to the point $\nu(u_1) = \nu(w_1)$ and see that $M_1$
is joined to either $\wt{N}$ or a
component $M_2$ dominating $\nu(\overline{w_1u_2})$ over $\nu(u_1)$.
If it is the former case, we again get a circuit in the dual graph of
$\wt{Y}_0$.
We may continue this argument and obtain that $\wt{N}$ is joined to
$M_1$ over $\nu(u_0)$, $M_1$ is joined to $M_2$ over $\nu(u_1)$ and so
on; finally, we have $M_{n-1}$ is joined to $M_n$ over $\nu(u_{n-1})$,
where $M_n\subset \wt{Y}_0$ is a component dominating
$\nu(\overline{w_{n-1} u_n})$. As mentioned before, there is no curve
$\overline{w_n u}\subset \wt{\Gamma}$. So $M_n$ is joined to $\wt{N}$
over $\nu(u_n) = \nu(w_n)$. Once again, we obtain a circuit in the dual
graph of $\wt{Y}_0$. Contradiction.

\figref{F5} illustrates our argument. Here \(\wt{N}\) passes through
\(u_0\) and \(w_i\). Then there will be a loop between \(\nu(u_0) =
\nu(w_0)\) and \(\nu(u_i) = \nu(w_i)\) on \(\wt{Y}_0\) and
consequently, \(p_a(\wt{Y}_0) > 0\). This is a contradiction.
\end{proof}


\begin{figure}[ht]
\centering
\begin{pspicture}(-10,-240)(430,20)

\rput[lt]{0}(0,0){
\begin{pspicture}(0,-240)(160,0)
\parametricplot{-10}{10}{80 8 t mul add
-10 t t mul 50 sub 150 t t mul sub mul 300 div add}
\parametricplot{-9}{9}{80 8 t mul add
-70 t t mul 50 sub 150 t t mul sub mul 300 div add}
\parametricplot{-9}{9}{80 8 t mul add
-100 t t mul 50 sub 150 t t mul sub mul 300 div add}
\parametricplot{-9}{9}{80 8 t mul add
-160 t t mul 50 sub 150 t t mul sub mul 300 div add}
\parametricplot{-10}{10}{80 8 t mul add
-220 t t mul 50 sub 150 t t mul sub mul 300 div add}

\uput[90](0,-2.33333){\small \(\wt{S}_\alpha\)}

\psline(0,-2.33333)(0,-212.33333)
\psline(160,-2.33333)(160,-212.333333)

\uput[90](40,-20.41666667){\small \(F_a\)}
\uput[100](120,-20.41666667){\small \(F_b\)}

\psline(40,-20.41666667)(40,-230.41666667)
\psline(120,-20.41666667)(120,-230.41666667)

\parametricplot{-0.5}{8}{80 t 8 mul sub
-125.41666667 t 15 mul add 25 t t mul sub t mul 9 div add}

\parametricplot[linestyle=dotted]{-2.5}{-0.5}{80 t 8 mul sub
-125.41666667 t 15 mul add 25 t t mul sub t mul 9 div add}

\parametricplot{-8}{-2.5}{80 t 8 mul sub
-125.41666667 t 15 mul add 25 t t mul sub t mul 9 div add}

\uput[45](55,-70){\small \(\wt{N}\)}

\pscircle*(40,-50.41666667){3}
\uput[225](40,-50.41666667){\small \(u_0\)}

\pscircle*(40,-80.41666667){3}
\uput[225](40,-80.41666667){\small \(u_1\)}

\pscircle*(40,-110.41666667){3}
\uput[225](40,-110.41666667){\small \(u_2\)}

\pscircle*(40,-170.41666667){3}
\uput[225](40,-170.41666667){\small \(u_i\)}

\pscircle*(120,-80.41666667){3}
\uput[315](120,-80.41666667){\small \(w_0\)}

\pscircle*(120,-110.41666667){3}
\uput[315](120,-110.41666667){\small \(w_1\)}

\pscircle*(120,-170.41666667){3}
\uput[315](120,-170.41666667){\small \(w_{i-1}\)}

\uput[225](110,-190){\small \(\wt{N}\)}

\pscircle*(120,-200.41666667){3}
\uput[350](120,-200.41666667){\small \(w_{i}\)}

\end{pspicture}
}


\psline{<->}(190,-110)(220,-110)

\rput[lt]{0}(240,0){
\begin{pspicture}(0,-210)(180,0)

\uput[90](20,-2.33333){\small \(\wt{N}\)}

\psline(20,-2.33333)(20,-212.333333)

\pscircle*(20,-30){3}
\uput[225](20,-30){\small \(\nu(u_0)\)}

\psline(0,-25)(120,-55)
\pscircle*(100,-50){3}
\uput[20](100,-50){\small \(\nu(u_1)\)}

\psline(105,-40)(80,-90)
\pscircle*(85,-80){3}
\uput[0](85,-80){\small \(\nu(u_2)\)}

\uput[90](60,-40){\small \(\overline{u_1 w_0}\)}

\pscircle*(20,-180){3}
\psline(0,-185)(120,-155)
\uput[135](20,-180){\small \(\nu(u_i)\)}

\psline(105,-170)(80,-120)
\pscircle*(85,-130){3}
\uput[0](85,-130){\small \(\nu(u_{i-2})\)}

\pscircle*(85,-95){2}
\pscircle*(85,-105){2}
\pscircle*(85,-115){2}

\pscircle*(100,-160){3}
\uput[340](100,-160){\small \(\nu(u_{i-1})\)}

\uput[0](130,-105){\small \(p_a(\wt{Y}_0) > 0\)}
\uput[300](60,-170){\small \(\overline{u_i w_{i-1}}\)}

\uput[270](50,-220){\small \(\wt{Y}_0\)}
\end{pspicture}
}

\end{pspicture}

\caption{\propref{prop5}}\label{F5}
\end{figure}

The difference between the points \(r_a, r_b\) and the other points
\(u_i, w_i\) lies in that at \(\nu(u_i) = \nu(w_i)\ne p_\alpha\),
\(X^{(\alpha)}\) is locally given by \(xy = t^\alpha\) so
\lemref{lem1} applies at \(\nu(u_i)\),
while \(X^{(\alpha)}\) has a rational double
point at \(p_\alpha = \nu(r_a) = \nu(r_b)\) and hence \lemref{lem1}
does not apply at \(\nu(r_a)\).

It is obvious that any two maximal S-chains are disjoint from each
other. Combining this with \propref{prop5}, we see that there is only
one maximal S-chain, i.e., the points in $\wt{\Gamma}\cap (F_a\cup
F_b)$ form an S-chain in a certain order. We can arrange the points in
$\wt{\Gamma}\cap (F_a\cup F_b)$ in the following way:
\begin{equation}\label{s4:e2}
\begin{split}
u_{-k}\xrightarrow{\varepsilon} w_{-k} & \xrightarrow{\phi} u_{-k+1}
\xrightarrow{\varepsilon} w_{-k+1}\xrightarrow{\phi} ...
\xrightarrow{\phi} u_0\xrightarrow{\varepsilon} w_0\\
&\xrightarrow{\phi} u_1\xrightarrow{\varepsilon} w_1 \xrightarrow{\phi} ...
\xrightarrow{\phi} u_l\xrightarrow{\varepsilon} w_l,
\end{split}
\end{equation}
where $u_0 = r_a$, $w_0 = r_b$ and $k, l\ge 0$.

\begin{prop}\label{prop6}
Let $\mu_i$ be the multiplicity of the curve $\overline{w_i u_{i+1}}$
in $\wt{\Gamma}$ for $-k\le i \le l-1$. Then
\begin{enumerate}
\item[A1.] $\mu_{-k}, \mu_{-k+1}, ..., \mu_0, \mu_1, ..., \mu_{l-1}$ satisfy
\begin{equation}\label{s4:e3}
1\le \mu_{-k}\le \mu_{-k+1}\le ...\le \mu_{-1} \text{ and }
\mu_0 \ge \mu_1 \ge ...\ge \mu_{l-1}\ge 1;
\end{equation}
\item[A2.] $Y^{(\alpha)}$ is totally separated at
$p_\alpha = \nu(u_0) = \nu(w_0)$ and hence
\begin{equation}\label{s4:e4}
|\mu_{-1} - \mu_0|\le 1;
\end{equation}
\item[A3.] if $N$ meets $\nu(\overline{w_i u_{i+1}})$ at a point
$s\ne\nu(w_i), \nu(u_{i+1})$, $Y^{(\alpha)}$ is totally separated at
$s$.
\end{enumerate}
\end{prop}

\begin{proof}
By \eqref{s4:e-1}, we have
\begin{equation}\label{s4:e5}
(\wt{\Gamma}\cdot F_a)_{u_i} = (\wt{\Gamma}\cdot F_b)_{w_i}
\end{equation}
for $-k\le i\le l$. So \eqref{s4:e3} is equivalent to the
statement that $\wt{N}$ meets $F_a$ only at the points
$u_{-k}, u_{-k+1}, ..., u_{-1}, u_0$ and meets $F_b$ only at the
points $w_0, w_1, ..., w_{l-1}, w_l$. Obviously, $\wt{N}$ must pass
through $u_{-k}$ since there is no curve $\overline{wu_{-k}}\subset
\wt{\Gamma}$. For the same reason, $w_l\in \wt{N}$.

Suppose that $w_{-i}\in \wt{N}$ for some $1\le i \le k$ and $i$ is
the largest number for this to hold. Applying \lemref{lem1} to
$\nu(w_{-i}) = \nu(u_{-i})$, we see that $\wt{N}$ is joined to a
component $M_1\subset \wt{Y}_0$ dominating 
\(\nu(\overline{w_{-i-1} u_{-i}})\) over
$\nu(w_{-i})$; continuing applying \lemref{lem1}, we see that $M_1$ is
joined to a component $M_2$ dominating
\(\nu(\overline{w_{-i-2} u_{-i-1}})\) over
$\nu(w_{-i-1})$, $M_2$ is joined to $M_3$ dominating
\(\nu(\overline{w_{-i-3}u_{-i -2}})\)
over $\nu(w_{-i-2})$ and so on. Finally, we have
$M_{k - i}$ dominating \(\nu(\overline{w_{-k} u_{-k+1}})\) 
is joined to $\wt{N}$ over
$\nu(w_{-k})$ and we obtain a circuit in the dual graph of
$\wt{Y}_0$. Contradiction. Therefore, \(w_{-k}, w_{-k+1}, ...,
w_{-1}\not\in \wt{N}\). Similarly, \(u_1, u_2, ..., u_l\not\in \wt{N}\).

If $Y^{(\alpha)}$ is not totally separated at $p_\alpha$, we have three
cases
\begin{enumerate}
\item a component $M_1\subset \wt{Y}_0$ dominating $\nu(\overline{w_0
u_1})$ is joined to a component $M_2\subset \wt{Y}_0$ dominating
$\nu(\overline{w_{-1} u_0})$ over $p_\alpha$;
\item a component $M_1\subset \wt{Y}_0$ dominating $\nu(\overline{w_0
u_1})$ is joined to $\wt{N}$ over $p_\alpha$;
\item a component $M_2\subset \wt{Y}_0$ dominating
$\nu(\overline{w_{-1} u_0})$ is joined to $\wt{N}$ over $p_\alpha$.
\end{enumerate}
In either of these cases, we can argue in the same way as before to
show that there is a circuit in the dual graph of $\wt{Y}_0$.
Therefore, $Y^{(\alpha)}$ is totally separated at
$p_\alpha$. As a consequence, by \coref{cor1}
$\wt{N}$ can be neither tangent to $F_a$ at $u_0$ nor tangent to
$F_b$ at $w_0$. So if $\wt{N}$ meets $F_a$
and $F_b$ at $u_0$ and $w_0$, it must meet $F_a$ and $F_b$
transversely at these points. Combining this with the fact that
$(\wt{\Gamma}\cdot F_a)_{u_0} = (\wt{\Gamma}\cdot F_b)_{w_0}$, we
obtain \eqref{s4:e4}.

Finally for (A3), if $Y^{(\alpha)}$ is not totally separated at $s =
N\cap \nu(\overline{w_i u_{i+1}})$, then $\wt{N}$ will be joined to a
component $M\subset \wt{Y}_0$ dominating $\nu(\overline{w_i u_{i+1}})$
over $s$. Again, we may use the same argument as before to show that
there is a circuit in the dual graph of $\wt{Y}_0$.
\end{proof}

Since $\wt{\Gamma} = \left(\cup_{i=-k}^{l-1} \mu_i \overline{w_i
u_{i+1}}\right) \cup \wt{N}$,
\begin{equation}\label{s4:e6}
\sum_{i=-\infty}^{\infty} \mu_i + \mu = m
\end{equation}
where $\mu$ is defined in \eqref{s4:e0} and
we let $\mu_i = 0$ if $i < -k$ or $i \ge l$. 
It follows from \eqref{s4:e5} that
\begin{equation}\label{s4:e7}
\mu_{i} - \mu_{i-1} = (\wt{N}\cdot F_a)_{u_i}
\end{equation}
for $i\le -1$ and
\begin{equation}\label{s4:e8}
\mu_j - \mu_{j+1} = (\wt{N}\cdot F_b)_{w_{j+1}}
\end{equation}
for $j\ge 0$. And since $\wt{N}$ meets $F_a$ and $F_b$ transversely at
$u_0$ and $w_0$ if it meets the curves at these points, we have
\begin{equation}\label{s4:e9}
\mu_0 \le \mu\le \mu_0 + 1\text{ and } \mu_{-1}\le \mu\le \mu_{-1} + 1
\end{equation}
where $\mu = \mu_0 + 1$ iff $w_0\in \wt{N}$ and $\mu = \mu_{-1} + 1$ iff
$u_0\in \wt{N}$. Hence
\begin{equation}\label{s4:e10}
(\wt{\Gamma}\cdot F_a)_{u_0} = (\wt{\Gamma}\cdot F_b)_{w_0} = \mu.
\end{equation}

Now we are ready to estimate the total $\delta$-invariant
$\delta(Y_t^{(\alpha)}, \Gamma)$ of $Y_t^{(\alpha)}$ in
the neighborhood of $\Gamma$. First, in the neighborhood of the
rational double point $p_\alpha$ where $Y^{(\alpha)}$
is totally separated by \propref{prop6}, we may apply \coref{cor2} to
conclude (noticing \eqref{s4:e10})
\begin{equation}\label{s4:e11}
\delta(Y_t^{(\alpha)}, p_\alpha) \ge \mu^2.
\end{equation}
Second, in the neighborhood of each point $s = N\cap \nu(\overline{w_i
u_{i+1}})$ with $s\not\in\{\nu(w_i), \nu(u_{i+1})\}$,
$Y^{(\alpha)}$ is totally separated by
\propref{prop6} and hence \lemref{lem2} can be applied (see also
\remref{rem3}). It follows that
\begin{equation}\label{s4:e12}
\delta(Y_t^{(\alpha)}, s) \ge \mu_i.
\end{equation}
Let $s_i = (N\cap \nu(\overline{w_i u_{i+1}}))\backslash
\{\nu(w_i), \nu(u_{i+1})\}$. Obviously, $s_i = \emptyset$ if either
$w_i\in \wt{N}$ or $u_{i+1}\in \wt{N}$. By \eqref{s4:e9},
$s_0 = \emptyset$ iff $\mu = \mu_0 + 1$. Therefore,
\begin{equation}\label{s4:e13}
\delta(Y_t^{(\alpha)}, s_0) \ge (\mu_0 + 1 - \mu)\mu_0
\end{equation}
by \eqref{s4:e12},
where we let $\delta(Y_t^{(\alpha)}, s_i) = 0$ if $s_i = \emptyset$.
Similarly,
\begin{equation}\label{s4:e13.5}
\delta(Y_t^{(\alpha)}, s_{-1}) \ge (\mu_{-1} + 1 - \mu)\mu_{-1}.
\end{equation}

Let $0\le a_0 < a_1 < a_2 < ... < a_n < ...$ 
be the sequence of integers such that
\begin{equation}\label{s4:e14}
\begin{split}
\mu_0 = ... = \mu_{a_0} & > \mu_{a_0+1} = \mu_{a_0+2} = ... =
\mu_{a_1}\\
& > \mu_{a_1+1} = \mu_{a_1+2} = ... = \mu_{a_2}\\
& > ... > \mu_{a_{n-1} + 1} = \mu_{a_{n-1} + 2} = ... = \mu_{a_n} > ... .
\end{split}
\end{equation}
Notice that for $i > 0$, $s_i = \emptyset$ iff $\mu_{i-1} \ne \mu_i$
by \eqref{s4:e8}. Therefore,
\begin{equation}\label{s4:e15}
\sum_{i > 0} \delta(Y_t^{(\alpha)}, s_i) \ge a_0\mu_0
+ \sum_{i > 0} (a_i - a_{i-1} - 1)\mu_{a_i}
\end{equation}
by \eqref{s4:e12}. Notice that
\begin{equation}\label{s4:e16}
\sum_{i\ge 0} \mu_i = (a_0 + 1)\mu_0 +\sum_{i>0} (a_i -
a_{i-1}) \mu_{a_i}.
\end{equation}
By \eqref{s4:e15} and \eqref{s4:e16},
\begin{equation}\label{s4:e17}
\begin{split}
\sum_{i > 0} \delta(Y_t^{(\alpha)}, s_i) - \sum_{i\ge 0} \mu_i &=
-\sum_{i\ge 0} \mu_{a_i}\\
&\ge -\left(\mu_0 + (\mu_0 - 1) + ... + 2 + 1\right)\\
&= -\frac{\mu_0(\mu_0 + 1)}{2}.
\end{split}
\end{equation}
By the same argument, we have
\begin{equation}\label{s4:e18}
\sum_{i < -1} \delta(Y_t^{(\alpha)}, s_i) - \sum_{i < 0} \mu_i
\ge -\frac{\mu_{-1}(\mu_{-1} + 1)}{2}.
\end{equation}
Putting \eqref{s4:e6}, \eqref{s4:e9}, \eqref{s4:e11}, \eqref{s4:e13},
\eqref{s4:e13.5}, \eqref{s4:e17} and \eqref{s4:e18} altogether, we obtain
\begin{equation}\label{s4:e19}
\begin{split}
\delta(Y_t^{(\alpha)}, \Gamma) &\ge \delta(Y_t^{(\alpha)}, p_\alpha) +
\delta(Y_t^{(\alpha)}, s_0) + \delta(Y_t^{(\alpha)}, s_{-1})\\
&\quad + \sum_{i > 0} \delta(Y_t^{(\alpha)}, s_i)
+ \sum_{i < -1} \delta(Y_t^{(\alpha)}, s_i)\\
&\ge m + \frac{1}{2} (\mu-\mu_0)^2 + \frac{1}{2} (\mu-\mu_{-1})^2\\
&\quad - \frac{1}{2} (\mu-\mu_0) - \frac{1}{2} (\mu-\mu_{-1}) = m.
\end{split}
\end{equation}
This finishes the proof of \eqref{s2:e10} and hence the first part of
\propref{prop4}.

It remains to find out what happens if $\delta(Y_t^{(\alpha)},
\Gamma) = m$.

\begin{prop}\label{prop7}
Suppose that $\delta(Y_t^{(\alpha)}, \Gamma) = m$. Then
\begin{enumerate}
\item[A1.] all the singularities of $Y_t^{(\alpha)}$ in the neighborhood of
$\Gamma$ actually lie in the neighborhoods of the points $p_\alpha$ and
$s_i$;
\item[A2.] the equality holds in \eqref{s4:e11};
\item[A3.] the equality holds in \eqref{s4:e12} for each $s = N\cap
\nu(\overline{w_i u_{i+1}})$ with $s\not\in\{\nu(w_i), \nu(u_{i+1})\}$;
\item[A4.] $\wt{N}$ meets $F_a$ and $F_b$ transversely at each
intersection, or equivalently,
\begin{equation}\label{s4:e20}
|\mu_i - \mu_{i+1}| \le 1
\end{equation}
for all $i$; in particular, $\mu_{-k} = \mu_{l-1} = 1$;
\item[A5.] for $-k\le i\le l-1$,
each component of $\wt{Y}_0$ that dominates $\nu(\overline{w_i
u_{i+1}})$ maps birationally to $\nu(\overline{w_i
u_{i+1}})$, i.e., there are no multiple
covers of $\nu(\overline{w_i u_{i+1}})$ on $\wt{Y}_0$.
\end{enumerate}
\end{prop}

\begin{figure}[ht]
\centering
\begin{pspicture}(0,-300)(470,20)

\rput[lt]{0}(0,0){
\begin{pspicture}(0,-300)(200,0)
\parametricplot{-10}{10}{100 10 t mul add
-10 t t mul 50 sub 150 t t mul sub mul 250 div add}
\parametricplot{-9}{9}{100 10 t mul add
-70 t t mul 50 sub 150 t t mul sub mul 250 div add}

\parametricplot{-9}{9}{100 10 t mul add
-130 t t mul 50 sub 150 t t mul sub mul 250 div add}
\parametricplot{-9}{9}{100 10 t mul add
-160 t t mul 50 sub 150 t t mul sub mul 250 div add}

\parametricplot{-9}{9}{100 10 t mul add
-220 t t mul 50 sub 150 t t mul sub mul 250 div add}

\parametricplot{-10}{10}{100 10 t mul add
-280 t t mul 50 sub 150 t t mul sub mul 250 div add}

\uput[90](0,0){\small \(\wt{S}_\alpha\)}

\psline(0,0)(0,-270)
\psline(200,0)(200,-270)

\uput[90](50,-22.5){\small \(F_a\)}
\uput[100](150,-22.5){\small \(F_b\)}

\psline(50,-22.5)(50,-292.5)
\psline(150,-22.5)(150,-292.5)

\pscircle*(50,-52.5){3}
\uput[225](50,-52.5){\small \(u_{-k}\)}

\pscircle*(50,-82.5){3}
\uput[225](50,-82.5){\small \(u_{-k+1}\)}

\uput[225](145,-250){\small \(\wt{N}\)}

\parametricplot{-7}{-3}{100 t 10 mul sub
-157.5 t 21 mul add 25 t t mul sub t mul 5 div add}

\parametricplot[linestyle=dotted]{-3}{-1.5}{100 t 10 mul sub
-157.5 t 21 mul add 25 t t mul sub t mul 5 div add}

\parametricplot{-1.5}{0.4}{100 t 10 mul sub
-157.5 t 21 mul add 25 t t mul sub t mul 5 div add}

\parametricplot[linestyle=dotted]{0.4}{2}{100 t 10 mul sub
-157.5 t 21 mul add 25 t t mul sub t mul 5 div add}

\parametricplot{2}{7}{100 t 10 mul sub
-157.5 t 21 mul add 25 t t mul sub t mul 5 div add}

\uput[45](65,-84){\small \(\wt{N}\)}

\pscircle*(50,-142.5){3}
\uput[225](50,-142.5){\small \(u_{0}\)}

\pscircle*(50,-172.5){3}
\uput[225](50,-172.5){\small \(u_{1}\)}

\pscircle*(50,-232.5){3}
\uput[225](50,-232.5){\small \(u_{l}\)}

\pscircle*(150,-82.5){3}
\uput[315](150,-82.5){\small \(w_{-k}\)}

\pscircle*(150,-142.5){3}
\uput[340](150,-142.5){\small \(w_{-1}\)}

\pscircle*(150,-172.5){3}
\uput[315](150,-172.5){\small \(w_0\)}

\pscircle*(150,-232.5){3}
\uput[315](150,-232.5){\small \(w_{l-1}\)}

\pscircle*(150,-262.5){3}
\uput[0](150,-262.5){\small \(w_{l}\)}
\end{pspicture}
}

\psline{<->}(220,-135)(250,-135)

\rput[lt]{0}(240,0){
\begin{pspicture}(0,-270)(230,0)

\uput[90](20,0){\small \(\wt{N}\)}

\psline(20,0)(20,-270)

\pscircle*(20,-30){3}
\uput[225](20,-30){\small \(u_{-k}\)}

\psline{<-}(43,-47)(90,0)
\uput[0](90,0){\small \(\overline{u_{-k+1} w_{-k}}\)}

\psline(10,-20)(70,-80)
\pscircle*(60,-70){3}

\psline{<-}(60,-80)(60,-100)
\uput[270](60,-100){\small \(u_{-k+1}\)}

\psline(50,-80)(110,-20)
\psline(130,-20)(190,-80)

\pscircle*(110,-50){2}
\pscircle*(120,-50){2}
\pscircle*(130,-50){2}

\pscircle*(180,-70){3}
\psline{<-}(180,-80)(180,-100)
\uput[270](180,-100){\small \(u_{-1}\)}

\psline(170,-80)(230,-20)
\uput[315](200,-50){\small \(\overline{u_0 w_{-1}}\)}

\pscircle*(220,-30){3}
\uput[135](220,-30){\small \(r_a\)}

\pscircle*(20,-240){3}
\uput[135](20,-240){\small \(w_l\)}

\psline{<-}(43,-223)(90,-270)
\uput[0](90,-270){\small \(\overline{u_l w_{l-1}}\)}

\psline(10,-250)(70,-190)
\pscircle*(60,-200){3}

\psline{->}(60,-170)(60,-190)
\uput[90](60,-170){\small \(w_{l-1}\)}

\psline(50,-190)(110,-250)
\psline(130,-250)(190,-190)

\pscircle*(110,-220){2}
\pscircle*(120,-220){2}
\pscircle*(130,-220){2}

\pscircle*(180,-200){3}

\psline{->}(180,-170)(180,-190)
\uput[90](180,-170){\small \(w_1\)}

\psline(170,-190)(230,-250)
\pscircle*(220,-240){3}
\uput[225](220,-240){\small \(r_b\)}

\uput[45](200,-220){\small \(\overline{u_1w_0}\)}

\uput[270](50,-280){\small \(\wt{Y}_0\)}
\end{pspicture}
}

\end{pspicture}

\caption{An admissible S-chain}\label{F6}
\end{figure}
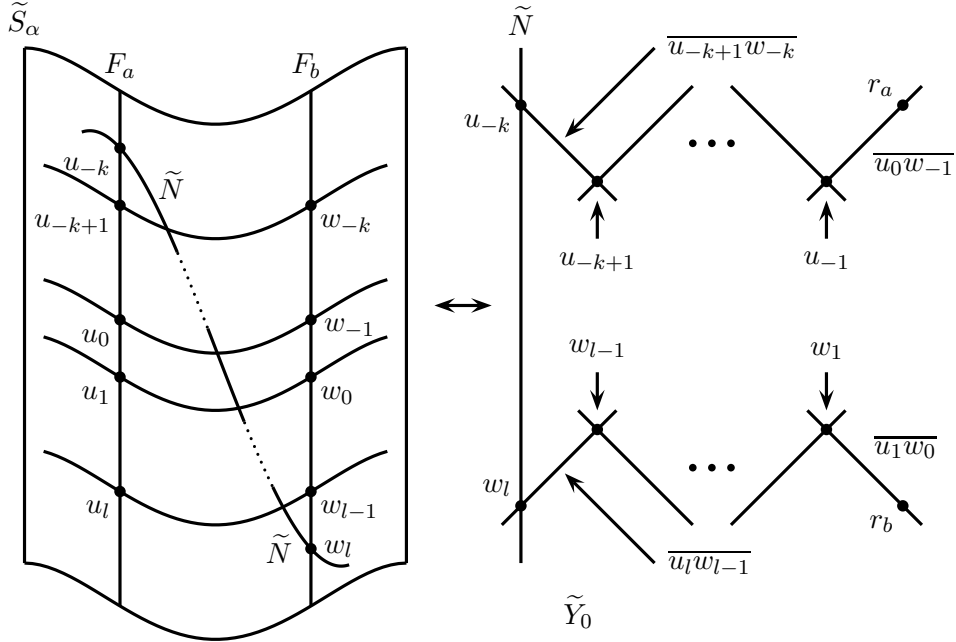

\begin{rem}\label{rem4}
In summary, the numerical relations among $\mu$ and $\mu_i$ are given
by \eqref{s4:e3}, \eqref{s4:e6}, \eqref{s4:e9} and \eqref{s4:e20}.
Those readers interested in the enumerative aspect of this
problem may have already noticed
that the number of such sequences $\{\mu, \mu_i\}$ can be
expressed in terms of partition numbers. As we already know, the
partition numbers have to pop up somewhere by the works of Yau-Zaslow
\cite{Y-Z} and Bryan-Leung \cite{B-L}. \figref{F6} shows the simplest
possible S-chain, corresponding to the case that
\(\mu_i = 1\) for \(-k\le i\le l-1\).
\end{rem}


\begin{proof}[Proof of \propref{prop7}]
Since $\delta(Y_t^{(\alpha)}, \Gamma) = m$, all the equalities in
\eqref{s4:e19} must hold. Then (A1), (A2) and (A3) follow immediately.

As for (A4), we notice that the equality in \eqref{s4:e17} has to hold.
So we must have $\mu_{a_0} = \mu_0$,
$\mu_{a_1} = \mu_0 - 1$, $\mu_{a_2} = \mu_0 - 2$ and so on, where
$\{a_i\}$ are defined by \eqref{s4:e14}. It follows immediately that
\eqref{s4:e20} holds for $i\ge 0$. Similarly, \eqref{s4:e20} holds for
$i < 0$. And by \eqref{s4:e7} and \eqref{s4:e8},
we see that $\wt{N}$ meets $F_a$ and
$F_b$ transversely everywhere.

Obviously, (A5) holds for $\nu(\overline{w_{l-1} u_l})$ and
$\nu(\overline{w_{-k} u_{-k+1}})$ since $\mu_{-k} = \mu_{l-1} = 1$.
Suppose that (A5) fails for some $\nu(\overline{w_i u_{i+1}})$ with $i
\ge 0$ and $i$ is the largest number with this property. Then there
exists a component $M\subset \wt{Y}_0$ dominating
$\nu(\overline{w_i u_{i+1}})$ with a map of degree at least $2$.
We claim that 
\begin{quote}
($*$) $M$ is joined to
at least two different components $M_1, M_2\subset \wt{Y}_0$ over
the point $\nu(u_{i+1})$, where $M_j = \wt{N}$ or $M_j$ dominates  
$\nu(\overline{w_{i+1} u_{i+2}})$ for $j=1,2$.
\end{quote}

If the map $M\to \nu(\overline{w_i u_{i+1}})$ is not totally ramified
over $\nu(u_{i+1})$, there are at least two distinct
points $x_1\ne x_2\in M$ such that $\pi(x_j) = \nu(u_{i+1})$ for
$j=1,2$ where $\pi: \wt{Y}\to Y^{(\alpha)}\subset X^{(\alpha)}$
is the map from $\wt{Y}$ to $Y^{(\alpha)}$. Then by \lemref{lem1}, the
branch of $M$ at $x_j$ is joined to a component $M_j$ over
the point $\nu(u_{i+1})$ for $j = 1, 2$, where $M_j = \wt{N}$ or
$\pi(M_j) = \nu(\overline{w_{i+1} u_{i+2}})$.
This justifies our claim ($*$) in the case that
$\pi: M\to \nu(\overline{w_i u_{i+1}})$ is not totally ramified over
$\nu(u_{i+1})$.

If $\pi: M\to \nu(\overline{w_i u_{i+1}})$ is totally ramified over
$\nu(u_{i+1})$, $\pi(M)$ meets $F_{p_{\alpha-1}}$ at $\nu(u_{i+1})$
with multiplicity at least $2$. Again by \lemref{lem1} (see also
\remref{rem2}), $M$ is joined to a union of components $\cup M_j$ over
$\nu(u_{i+1})$ such that $\pi(\cup M_j)\subset \nu(\overline{w_{i+1}
u_{i+2}}) \cup N$ and $\pi(\cup M_j)$ meets $F_{p_{\alpha-1}}$ at
$\nu(u_{i+1})$ with multiplicity at least $2$. Our assumption on $i$
implies that (A5) holds for $\nu(\overline{w_{i+1} u_{i+2}})$, i.e.,
every component of $\wt{Y}_0$ dominating $\nu(\overline{w_{i+1}
u_{i+2}})$ maps birationally to $\nu(\overline{w_{i+1}
u_{i+2}})$. And since $\wt{N}$
meets $F_b$ transversely at $w_{i+1}$ if $w_{i+1}\in \wt{N}$, we see
that $\cup M_j$ contains at least two different components dominating
either $N$ or $\nu(\overline{w_{i+1} u_{i+2}})$ and hence ($*$) follows.

Starting with ($*$), we may argue as before to show that each $M_j$ is
joined by a chain of components over
$\nu(\overline{w_{i+2} u_{i+3}}) \cup \nu(\overline{w_{i+3}
u_{i+4}})\cup ...\cup \nu(\overline{w_{l-1} u_{l}})$ to $\wt{N}$ for
$j = 1, 2$. And hence there is a circuit in the dual graph of
$\wt{Y}_0$. Contradiction. So (A5) holds for each $\nu(\overline{w_i
u_{i+1}})$ with $i\ge 0$. A similar argument shows that (A5) holds for each
$\nu(\overline{w_i u_{i+1}})$ with $i< 0$. 
\end{proof}

With \propref{prop7}, the second part of \propref{prop4} is almost
immediate. In the neighborhood of $p_\alpha$, $Y^{(\alpha)}$ consists
of $2\mu$ local irreducible components corresponding to $2\mu$ branches
of $\wt{Y}_0$ over $p_\alpha$. And since the equality holds in
\eqref{s4:e11}, $Y_t^{(\alpha)}$ has exactly $\mu^2$ nodes as
singularities in the neighborhood of $p_\alpha$ by \coref{cor2}.
In the neighborhood of a point $s = N\cap
\nu(\overline{w_i u_{i+1}})$ with $s\not\in\{\nu(w_i),
\nu(u_{i+1})\}$, $Y^{(\alpha)}$ consists of $\mu_i + 1$ local
irreducible components. And since the equality holds in
\eqref{s4:e12}, $Y_t^{(\alpha)}$ has exactly $\mu_i$ nodes as
singularities in the neighborhood of $s$ by \lemref{lem2}. Therefore,
$Y_t^{(\alpha)}$ is nodal if $\delta(Y_t^{(\alpha)}, \Gamma) = m$. This
finishes the proof of \propref{prop4}.

Although it is no longer necessary for our purpose, it will be interesting
to classify all the possible configurations for the stable reduction
$\wt{Y}_0$. Actually, this is not hard given everything we have done
so far. Next, we will give a description for $\wt{Y}_0$ without
justification and leave the readers to verify the details.

Let us contract some curves on $\wt{Y}_0$ to make $\wt{Y}\to Y$ into a
stable map. Remember that we start with the stable map $\wt{Y}\to
Y^{(\alpha)}$.

Among the components of $\wt{Y}_0$ that dominate $E$,
\begin{enumerate}
\item there is only one component $\wt{N}$ dominating $E$ with
a map of degree $\mu$ and the rest each map to $E$ birationally;
\item the map $\wt{N}\to \wt{E}$ is unramified over $a$ and $b$, where
$\wt{E}$ is the normalization of $E$ and $a, b\in \wt{E}$
are the two points over the node $p\in E$;
\item two components $M_1$ and $M_2$ only meet at a point $x$ over the
node $p$; in addition, $M_1\cup M_2$ maps biholomorphically to $E$ locally
at $x$, i.e., the two branches of $M_1$ and $M_2$ at $x$ must map to
different branches of $E$ at $p$; using the terminology of \cite{B-L}, we
say that there is a ``branch jump'' whenever two components meet;
\item for each $x\in \wt{N}$ over $p$, there is a chain of curves
$\cup M_i$ attached to $\wt{N}$ at $p$ with each $M_i$ dominating $E$;
and each component $M\ne \wt{N}$ dominating $E$
lies on one of these $2\mu$ chains;
\item let $\lambda(x)$ be the length of the chain of curves attached
to the point $x\in \wt{N}$ over $p$; obviously,
\begin{equation}\label{s4:e21}
\sum \lambda(x) + \mu = m
\end{equation}
where we sum over all the $2\mu$ points $x\in \wt{N}$ that map to $p$;
\item for any two points $x_1\ne x_2\in \wt{N}$ over $a$,
$\lambda(x_1) \ne \lambda(x_2)$; similarly, for any two points $y_1\ne
y_2\in \wt{N}$ over $b$, $\lambda(y_1) \ne \lambda(y_2)$;
\item $\wt{N}$ meets $\wt{C}$ at a point over $q = C\cap E$, where
$\wt{C}$ is the component of $\wt{Y}_0$ dominating $C$.
\end{enumerate}

Let \(x_1, x_2, ..., x_\mu\) be the points of \(\wt{N}\) over \(a\)
and \(y_1, y_2, ..., y_\mu\) be the points of \(\wt{N}\) over
\(b\). Then \(\{x_i\}\) map to the points among
\(u_{-k}, u_{-k+1}, ..., u_l\) and \(\{y_i\}\) map to the
points among \(w_{-k}, w_{-k+1}, ..., w_l\).
Let \(\lambda_i = \lambda(x_i)\) and \(\lambda_{-i} = \lambda(y_i)\) and
we order \(x_i\) and \(y_i\) such that
\begin{equation}\label{s4:e22}
\lambda_1 > \lambda_2 > ... > \lambda_\mu \ge 0
\end{equation}
and
\begin{equation}\label{s4:e23}
\lambda_{-1} > \lambda_{-2} > ... > \lambda_{-\mu} \ge 0
\end{equation}
where \(\lambda_\mu = 0\) if and only if \(x_\mu\) maps to \(u_0 =
r_a\) and \(\lambda_{-\mu} = 0\) if and only if \(y_\mu\) maps to
\(w_0 = r_b\). Under these notations, we may rewrite \eqref{s4:e21} as
\begin{equation}\label{s4:e24}
\sum_{i=-\mu}^\mu \lambda_i + \mu = m
\end{equation}
where we let \(\lambda_0 = 0\). Later in \appendixref{apdx2}
when we count the number of rational curves on a K3 surface,
we are basically
counting the number of the sequences \(\{\mu, \lambda_i\}\)
satisfying \eqref{s4:e22}, \eqref{s4:e23} and
\eqref{s4:e24}. \figref{F7} shows the configuration of
\(\wt{Y}_0\).
Also see \figref{F6} for the
simplest possible configuration of \(\wt{Y}_0\), corresponding to
the case that \(\mu = 1\).


\begin{figure}[ht]
\centering
\begin{pspicture}(0,-420)(240,20)

\uput[90](20,0){\small \(\wt{N}\)}

\psline(20,0)(20,-420)

\uput[225](20,-10){\small \(x_1\)}

\psline(15,-5)(40,-30)
\psline(30,-30)(55,-5)

\psline(45,-5)(70,-30)
\psline(60,-30)(85,-5)

\psline(75,-5)(100,-30)
\psline(90,-30)(115,-5)

\psline(105,-5)(130,-30)
\psline(120,-30)(145,-5)

\psline(195,-5)(220,-30)
\psline(210,-30)(235,-5)

\pscircle*(160,-17.5){2}
\pscircle*(170,-17.5){2}
\pscircle*(180,-17.5){2}

\psline{<-}(20,-45)(107.5,-45)
\psline{->}(147.5,-45)(235,-45)
\rput{0}(127.5,-45){\small \(\lambda_1\)}

\uput[225](20,-60){\small \(x_2\)}

\psline(15,-55)(40,-80)
\psline(30,-80)(55,-55)

\psline(45,-55)(70,-80)
\psline(60,-80)(85,-55)

\psline(75,-55)(100,-80)
\psline(90,-80)(115,-55)

\psline(165,-55)(190,-80)
\psline(180,-80)(205,-55)

\pscircle*(130,-67.5){2}
\pscircle*(140,-67.5){2}
\pscircle*(150,-67.5){2}

\psline{<-}(20,-95)(92.5,-95)
\psline{->}(132.5,-95)(205,-95)
\rput{0}(112.5,-95){\small \(\lambda_2\)}

\pscircle*(80,-120){2}
\pscircle*(80,-130){2}
\pscircle*(80,-140){2}

\uput[225](20,-160){\small \(x_\mu\)}

\psline(15,-155)(40,-180)
\psline(30,-180)(55,-155)

\psline(105,-155)(130,-180)
\psline(120,-180)(145,-155)

\pscircle*(70,-167.5){2}
\pscircle*(80,-167.5){2}
\pscircle*(90,-167.5){2}

\psline{<-}(20,-195)(62.5,-195)
\psline{->}(102.5,-195)(145,-195)
\rput{0}(82.5,-195){\small \(\lambda_\mu\)}

\rput[lt]{0}(-5,-210){
\begin{pspicture}(0,-210)(240,0)
\uput[225](20,-10){\small \(y_1\)}

\psline(15,-5)(40,-30)
\psline(30,-30)(55,-5)

\psline(45,-5)(70,-30)
\psline(60,-30)(85,-5)

\psline(75,-5)(100,-30)
\psline(90,-30)(115,-5)

\psline(105,-5)(130,-30)
\psline(120,-30)(145,-5)

\psline(195,-5)(220,-30)
\psline(210,-30)(235,-5)

\pscircle*(160,-17.5){2}
\pscircle*(170,-17.5){2}
\pscircle*(180,-17.5){2}

\psline{<-}(20,-45)(107.5,-45)
\psline{->}(147.5,-45)(235,-45)
\rput{0}(127.5,-45){\small \(\lambda_{-1}\)}

\uput[225](20,-60){\small \(y_2\)}

\psline(15,-55)(40,-80)
\psline(30,-80)(55,-55)

\psline(45,-55)(70,-80)
\psline(60,-80)(85,-55)

\psline(75,-55)(100,-80)
\psline(90,-80)(115,-55)

\psline(165,-55)(190,-80)
\psline(180,-80)(205,-55)

\pscircle*(130,-67.5){2}
\pscircle*(140,-67.5){2}
\pscircle*(150,-67.5){2}

\psline{<-}(20,-95)(92.5,-95)
\psline{->}(132.5,-95)(205,-95)
\rput{0}(112.5,-95){\small \(\lambda_{-2}\)}

\pscircle*(80,-120){2}
\pscircle*(80,-130){2}
\pscircle*(80,-140){2}

\uput[225](20,-160){\small \(y_\mu\)}

\psline(15,-155)(40,-180)
\psline(30,-180)(55,-155)

\psline(105,-155)(130,-180)
\psline(120,-180)(145,-155)

\pscircle*(70,-167.5){2}
\pscircle*(80,-167.5){2}
\pscircle*(90,-167.5){2}

\psline{<-}(20,-195)(62.5,-195)
\psline{->}(102.5,-195)(145,-195)
\rput{0}(82.5,-195){\small \(\lambda_{-\mu}\)}
\end{pspicture}
}
\end{pspicture}

\caption{\(\wt{Y}_0\)}\label{F7}
\end{figure}

\begin{figure}[ht]
\centering
\begin{pspicture}(-20,-260)(430,30)

\uput[90](20,0){\small \(\lambda_1\)}
\uput[90](60,0){\small \(\lambda_2\)}
\uput[90](100,0){\small \(\lambda_3\)}
\uput[90](140,0){\small \(\lambda_4\)}
\uput[90](180,0){\small \(\lambda_5\)}

\uput[180](0,-20){\small \(\mu_{-1}\)}
\uput[180](0,-60){\small \(\mu_{-2}\)}
\uput[180](0,-100){\small \(\mu_{-3}\)}
\uput[180](0,-140){\small \(\mu_{-4}\)}
\uput[180](0,-180){\small \(\mu_{-5}\)}
\uput[180](0,-220){\small \(\mu_{-6}\)}

\psline(0,0)(200,0)
\psline(0,-40)(200,-40)
\psline(0,-80)(160,-80)
\psline(0,-120)(160,-120)
\psline(0,-160)(120,-160)
\psline(0,-200)(80,-200)
\psline(0,-240)(40,-240)

\psline(0,0)(0,-240)
\psline(40,0)(40,-240)
\psline(80,0)(80,-200)
\psline(120,0)(120,-160)
\psline(160,0)(160,-120)
\psline(200,0)(200,-40)

\rput[lt]{0}(230,0){
\begin{pspicture}(0,-260)(200,0)
\uput[90](20,0){\small \(\lambda_{-1}\)}
\uput[90](60,0){\small \(\lambda_{-2}\)}
\uput[90](100,0){\small \(\lambda_{-3}\)}
\uput[90](140,0){\small \(\lambda_{-4}\)}
\uput[90](180,0){\small \(\lambda_{-5}\)}

\uput[180](0,-20){\small \(\mu_{0}\)}
\uput[180](0,-60){\small \(\mu_{1}\)}
\uput[180](0,-100){\small \(\mu_{2}\)}
\uput[180](0,-140){\small \(\mu_{3}\)}
\uput[180](0,-180){\small \(\mu_{4}\)}
\uput[180](0,-220){\small \(\mu_{5}\)}

\psline(0,0)(200,0)
\psline(0,-40)(200,-40)
\psline(0,-80)(200,-80)
\psline(0,-120)(160,-120)
\psline(0,-160)(120,-160)
\psline(0,-200)(80,-200)
\psline(0,-240)(40,-240)

\psline(0,0)(0,-240)
\psline(40,0)(40,-240)
\psline(80,0)(80,-200)
\psline(120,0)(120,-160)
\psline(160,0)(160,-120)
\psline(200,0)(200,-80)
\end{pspicture}
}
\end{pspicture}

\caption{Relation between \(\{\lambda_i\}\) and \(\{\mu_j\}\)}\label{F8}
\end{figure}
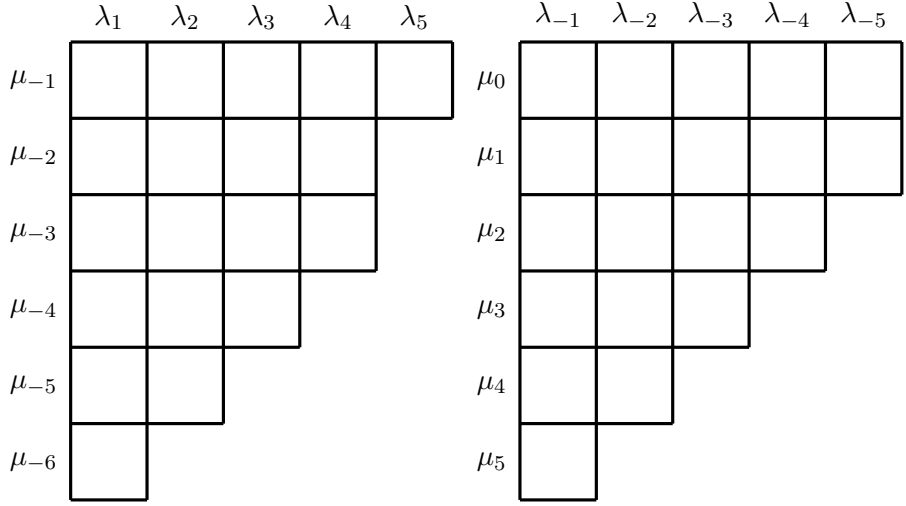

It is also worthwhile to point out that \(\{\mu, \lambda_i\}\)
are uniquely determined by \(\{\mu, \mu_j\}\) and vice
versa. Actually, we can describe their relation explicitly as
follows: the Young tableau of \((\lambda_1,
\lambda_2,...,\lambda_\mu)\) is dual to the Young tableau of
\((\mu_{-1}, \mu_{-2}, ..., \mu_{-k})\) and the Young tableau of
\((\lambda_{-1},
\lambda_{-2},...,\lambda_{-\mu})\) is dual to the Young tableau of
\((\mu_0, \mu_1, ..., \mu_{l-1})\) (see \figref{F8}).


\section{Proofs of \propref{prop2} and \ref{prop3}}\label{s5}

\subsection{Proof of \propref{prop2}}\label{s5:1}

Although the proposition says that we make a base change of a
one-size-fits-all degree $\alpha$ at the very beginning, in practice
we have no idea of what values $\alpha$ should take before we start to blow
up $X$ and $Y$. So our proof goes as follows: we start with an
$\alpha$ for which the proposition might fail, then we make a sequence
of base changes depending on where it fails and finally we will obtain
an $\alpha$ such that everything in the proposition holds.

Suppose that the proposition holds for $Y^{(n)}$
and $Y_0^{(n)}$ contains $E_n$ with multiplicity $\mu$. So we start
with $n = 0$ and $\mu = m$ and we will show that eventually either $n
= \alpha$ or $\mu = 0, 1$.

Suppose that $\mu \ge 2$. Pick an arbitrary smooth point $b\ne q_n \in
E_n$ of $E_n$. Locally at $b$, the curve $Y_0^{(n)}$ is given by
$z^\mu = 0$ in $\Delta_{wz}^2$. With a suitable choice of the
coordinate $z$, the family $Y^{(n)}$ is locally given by
\begin{equation}\label{s5:e1}
z^\mu + t^{a_1} f_1(t, w) z^{\mu - 2} + t^{a_2} f_2(t, w) z^{\mu - 3}
+ ... + t^{a_{\mu - 1}} f_{\mu - 1}(t, w) = 0
\end{equation}
in $\Delta_{wzt}^3$, where $a_i > 0$ and $f_i(0, 0) \ne 0$ for
$i=1,2,...,\mu-1$. Let
\begin{equation}\label{s5:e2}
\beta = \min_{1\le i \le \mu - 1}\frac{a_i}{i + 1}.
\end{equation}
A base change might be needed in order to make $\beta$ into a positive
integer and we have to modify the sequence \eqref{s2:e9} after a base
change. But the bottom line is that $\mu$ does not change in the
process. So let us assume that $\beta$ is a positive integer.

If $\beta > 1$, the local equation \eqref{s5:e1} shows that
$M = Y^{(n+1)}\cap S_{n+1}$ contains a section of $S_{n+1}\to E_n$ with
multiplicity $\mu$ and $E_n \not\subset M$.
And local computations of $Y^{(n)}$ at $p_n$ and $q_n$ show that
$M$ meets $E_n$ only at $q_n$ and it meets $E_n$ at $q_n$ transversely.
So $M = F_{q_n}\cup \mu G$, where $G\in |\CO_{S_{n+1}}(E_n)|$ and
$G\ne E_n$. If $G\ne
E_{n+1}$, we are done with the proof of \propref{prop2} since any
further blowups of $Y^{(n+1)}$ will only produce more $F_{q_i}$'s,
i.e., $Y^{(n+k)}\cap S_{n+k}$ will consist only of $F_{q_{n+k-1}}$ for
$k > 1$.

So we can apply this argument to every $1\le k \le \beta - 1$: either
$Y^{(n+k)}\cap S_{n+k} = F_{q_{n+k-1}}\cup \mu E_{n+k}$ for each $k$
or this fails for certain $k$ such that $Y^{(n+k)}\cap S_{n+k} =
F_{q_{n+k-1}} \cup \mu G$ with $G\ne E_{n+k}$, in which case we are done.

Let us assume that $Y^{(n+k)}\cap S_{n+k} = F_{q_{n+k-1}}\cup \mu
E_{n+k}$ for each $k=1,2,...,\beta-1$.

Due to our choice \eqref{s5:e2} of $\beta$,
$M = Y^{(n+\beta)}\cap S_{n+\beta}$
consists of at least two components, each of which dominates $E_{n+\beta -
1}$ with a map of degree strictly less than $\mu$, and
$E_{n+\beta-1}\not\subset M$. Again, local
computations of $Y^{(n+\beta-1)}$ at $p_{n+\beta-1}$ and
$q_{n+\beta-1}$ show that $M$ meets $E_{n+\beta-1}$ only at
$q_{n+\beta-1}$ and it meets $E_{n+\beta-1}$ at $q_{n+\beta-1}$
transversely. So $M\in |\CO_{S_{n+\beta}}(F_{q_{n+\beta-1}} + \mu
E_{n+\beta-1})|$. It remains to show that $F_{q_{n+\beta-1}}\subset M$ if
$n+\beta < \alpha$.

Let $\nu: \wt{E}_{n+\beta-1}\to E_{n+\beta-1}$ be the normalization of
$E_{n+\beta-1}$. It induces the normalization $\nu: \P^1\times
\wt{E}_{n+\beta-1} \to S_{n+\beta} \isom \P^1\times E_{n+\beta - 1}$ of
$S_{n+\beta}$. Let $a, b\in \wt{E}_{n+\beta-1}$ be the preimages of
$p_{n+\beta-1}$ and let $F_a$ and $F_b$ be the fibers over
$a$ and $b$.

Let $\phi_{ab}$ be the natural identification between $F_a$ and $F_b$
as defined in \ssecref{s3:1}. We can think of $S_{n+\beta}$ as
obtained from $\P^1\times \wt{E}_{n+\beta-1}$ by gluing $F_a$ and
$F_b$ via $\phi_{ab}$. Let $r_a\in F_a$ and $r_b\in F_b$ be the
preimages of the rational double point $p_{n+\beta}$ of
$X^{(n+\beta)}$. Obviously, $\phi_{ab}(r_a) = r_b$.

Let $\wt{M} = \nu^{-1}(M)$. If $\wt{M}$ meets $F_a$ at a point $s_a\ne
r_a$ with multiplicity $k$, the branches of $\wt{M}$ at $s_a$ will map
to the branches of $M$ lying on one of two surfaces of
$X_0^{(n+\beta)}$ at $\nu(s_a)$. Recall
that $X^{(n+\beta)}$ is locally given by $xy = t^{n+\beta}$ at
$\nu(s_a)$. So we can apply \lemref{lem1} to conclude that there exist
branches of $M$ lying on the other surface of $X_0^{(n+\beta)}$
at $\nu(s_a)$ and the branches on
both surfaces meet $F_{p_{n+\beta-1}}$ at $\nu(s_a)$ with the same
multiplicity $k$. Correspondingly, $\wt{M}$ must meet $F_b$ at
$s_b = \phi_{ab}(s_a)$ with the same multiplicity $k$. And since $\wt{M}\in
|F_{q_{n+\beta-1}} + \mu \wt{E}_{n+\beta-1}|$, we draw the conclusion
that if $\wt{M}$ meets $F_a$ at a point $s_a\ne r_a$ with multiplicity
$k$, it must meet $F_b$ at $s_b$ with the same multiplicity
$k$ and hence it must contain the curve $\overline{s_a s_b}$
with multiplicity $k$. Similarly, if $\wt{M}$ meets $F_b$ at
a point $s_b\ne r_b$ with multiplicity $k$,
$\wt{M}$ must meet $F_a$ at $s_a = \phi_{ba}(s_b)$ with the same
multiplicity $k$ and hence it must contain the curve $\overline{s_a s_b}$
with multiplicity $k$. Therefore, if we let $\wt{N}\subset \wt{M}$ be the
irreducible component of $\wt{M}$ with
$\wt{N}\in |F_{q_{n+\beta-1}} + \gamma \wt{E}_{n+\beta-1}|$ for some
$\gamma \le \mu$, we see that $\wt{N}$ meets $F_a$ and $F_b$ only at $r_a$
and $r_b$. But then $\overline{r_a r_b}\subset \wt{N}$ if $\gamma >
0$. Therefore, $\gamma = 0$ and $F_{q_{n+\beta-1}}\subset M$.

Since $M$ has at least two components which dominates $E_{n+\beta-1}$,
the multiplicity $\mu'$ of $E_{n+\beta}$ in $M$ is strictly less than
$\mu$. Now the proposition holds for $Y^{(n+\beta)}$ and
$Y_0^{(n+\beta)}$ contains $E_{n+\beta}$ with multiplicity $\mu' < \mu$.
We see that the value of $\mu$ has been reduced.

Finally, if $\mu = 0$, there is nothing left to do. If $\mu = 1$,
no further base changes are needed; we just have to verify that
$F_{q_{n+k - 1}} \subset M = Y^{(n+k)}\cap S_{n+k}$ for $1 \le k \le
\alpha - n - 1$, the argument for which goes exactly as before. This
completes the proof of \propref{prop2}.

\subsection{Proof of \propref{prop3}}

Suppose that $Y^{(\alpha)} \cap S_n$ contains a component $M\in
|\CO_{S_n}(E_{n-1})|$ with multiplicity $\mu > 0$ for some $1\le n \le
\alpha - 1$. Namely, $M$ is a wondering component. Let $u\in M$ be the
node of $M$, where $X^{(\alpha)}$ is locally given by $xy =
t^n$.

Let $\wt{Y}\to Y^{(\alpha)}$ be the stable reduction of
$Y^{(\alpha)}$ after normalization defined as before.
Let $G$ be the dual graph of the components of $\wt{Y}_0$ that
map to $M$ (including the curves contracted to a point on $M$)
and let us remove from $G$ all the vertices of degree $0$ or $1$ that
represent contractible curves. So $\deg([R]) \ge 2$ for any $[R]\in G$
representing a contractible curve $R$, where we let $[A]$ denote the
vertex of $G$ representing the component $A\subset \wt{Y}_0$ and
$\deg([A])$ denote the degree of $[A]$ in $G$.

Let $\wt{M}\subset \wt{Y}_0$ be a component of $\wt{Y}_0$ dominating
$M$ and let $\wt{u}\in \wt{M}$ be one of the points over the node
$u$. The branch of $\wt{M}$ at $\wt{u}$ maps to a branch of $M$ lying
on one of two surfaces of $X_0^{(\alpha)}$ at $u$. So by
\lemref{lem1}, the branch of $\wt{M}$ at $\wt{u}$ is joined by a chain
of contractible curves to a branch of $\wt{Y}_0$ that maps to the branch
of $M$ lying on the other surface of $X_0^{(\alpha)}$ at $u$. This is
to say that each $\wt{u}\in \wt{M}$ over $u$ corresponds to an edge of
$G$ from $[\wt{M}]$.
And since there are at least two
points of $\wt{M}$ mapping to $u$, we must have $\deg([\wt{M}]) \ge 2$
in $G$.

So every vertex of $G$ has degree at least $2$. This is impossible and
hence \propref{prop3} follows.

\appendix

\section{Deformations of K3 Surfaces}\label{apdx1}

Here we will give a proof for \lemref{lem0}.

Without the loss of generality, let us assume that \(D\) is very
ample; otherwise, we may simplely replace \(Y\) by \(nY\) and \(D\) by
\(nD\) for some \(n>>0\). Let \(g = \dim |\CO_X(Y)| = \dim
|\CO_S(D)|\) and \(X\) can be embedded to \(\P^g\times \Delta\) by the
complete linear series \(|\CO_X(Y)|\).
Let \(N_S\) be the normal bundle of \(S\) in \(\P^g\)
and we have the standard exact sequence
\begin{equation}\label{lem0:e2}
0\xrightarrow{} T_S\xrightarrow{} T_{\P^g} |_S \xrightarrow{} N_S \xrightarrow{} 0
\end{equation}
on \(S\). From \eqref{lem0:e2}, we have the exact sequence
\begin{equation}\label{lem0:e3}
H^0(N_S) \xrightarrow{} H^1(T_S) \xrightarrow{f} H^1(T_{\P^g} |_S).
\end{equation}
The embedding \(X\hookrightarrow \P^g\times \Delta\) gives an embedded
deformation of \(S\) in \(\P^g\). Therefore, the Kodaira-Spencer map
\(\ks: T_{\Delta,0}\to H^1(T_S)\) factors through
\(H^0(N_S)\). Consequently, \(\ks(\partial/\partial t)\) lies in the
kernel of the map \(f: H^1(T_S) \to H^1(T_{\P^g} |_S)\). Therefore, to
prove \eqref{lem0:e1}, it suffices to show that
\begin{equation}\label{lem0:e4}
\ker f = V.
\end{equation}

The Euler sequence
\begin{equation}\label{lem0:e5}
0\xrightarrow{} \CO_S \xrightarrow{} \left. \CO_{\P^g}(1)^{\oplus
(g+1)}\right|_S \xrightarrow{} T_{\P^g} |_S\xrightarrow{} 0
\end{equation}
yields an isomorphism from \(H^1(T_{\P^g} |_S)\) to \(H^2(\CO_S)\)
since
\begin{equation}\label{lem0:e6}
H^i(\CO_{\P^g}(1) |_S) = H^i(\CO_S(D)) = 0
\end{equation}
for \(i=1,2\) by Kodaira vanishing theorem. So we have
\begin{equation}\label{lem0:e7}
H^1(T_S) \xrightarrow{f} H^1(T_{\P^g} |_S) \xrightarrow{\sim} H^2(\CO_S).
\end{equation}
Let us consider the dual sequence of \eqref{lem0:e7}, i.e.,
\begin{equation}\label{lem0:e8}
\begin{CD}
H^1(T_S) @>f>> H^1(T_{\P^g} |_S) @>\sim>> H^2(\CO_S)\\
\times &&\times &&\times\\
H^1(\Omega_S) @<f^\vee<< H^1(\Omega_{\P^g} |_S) @<\sim<< H^0(\CO_S).\\
@VVV @VVV @VVV\\
\BC && \BC && \BC
\end{CD}
\end{equation}
Obviously, \eqref{lem0:e4} is equivalent to saying that the image of
the map \(f^\vee: H^1(\Omega_{\P^g} |_S)\to H^1(\Omega_S)\) is spanned by
\(c_1(D)\). So it suffices to prove that
\begin{equation}\label{lem0:e9}
\im f^\vee = \Span\{c_1(D)\}.
\end{equation}
From the dual Euler sequences
\begin{equation}\label{lem0:e10}
\begin{CD}
0 @>>> \Omega_{\P^g} @>>> \CO_{\P^g}(-1)^{\oplus
(g+1)} @>>> \CO_{\P^g} @>>> 0\\
&& @VVV @VVV @VVV\\
0 @>>> \Omega_{\P^g} |_S @>>> \left. \CO_{\P^g}(-1)^{\oplus
(g+1)}\right|_S @>>> \CO_S @>>> 0
\end{CD}
\end{equation}
we see that
\begin{equation}\label{lem0:e11}
\begin{CD}
H^0(\CO_{\P^g}) @>{\sim}>> H^1(\Omega_{\P^g})\\
@VV{\sim}V @VV{\sim}V\\
H^0(\CO_S) @>{\sim}>> H^1(\Omega_{\P^g} |_S).
\end{CD}
\end{equation}
It is a common knowledge that \(H^1(\Omega_{\P^g}) = \BC\) is generated by
\(c_1\left(\CO_{\P^g}(1)\right)\). It is not hard to see that the
image of \(c_1\left(\CO_{\P^g}(1)\right)\) under the map
\begin{equation}\label{lem0:e12}
H^1(\Omega_{\P^g}) \xrightarrow{\sim} H^1(\Omega_{\P^g} |_S)
\xrightarrow{f} H^1(\Omega_S)
\end{equation}
is \(c_1(D)\). This proves \eqref{lem0:e9} and
\eqref{lem0:e9} \(\Rightarrow\) \eqref{lem0:e4} \(\Rightarrow\)
\eqref{lem0:e1}.

For the second part of the lemma,
suppose that \(S\) is a K3 surface, \(D\) is an ample divisor on \(S\)
and \(S\) is embedded into \(\P^g\) by \(|\CO_S(nD)|\) for some \(n > 0\).
Observe that \eqref{lem0:e4} also implies that \(f\)
is a surjection and hence \(H^1(N_S) = 0\) by the
exact sequence
\begin{equation}\label{lem0:e13}
H^1(T_S) \xrightarrow{f} H^1(T_{\P^g}|_S) \xrightarrow{} H^1(N_S)
\xrightarrow{} H^2(T_S) = 0.
\end{equation}
So the embedded deformations of \(S\) in \(\P^g\) are
unobstructed. And since \(H^0(N_S)\) surjects onto \(V\), for each
\(v\in V\), there exists an embedded deformation of \(S\subset \P^g\) with
Kodaira-Spencer class \(v\), i.e., there exists a smooth family \(X\)
over \(\Delta\) and an embedding \(\varphi: X\hookrightarrow
\P^g\times \Delta\) such that \(\varphi(X_0) = S\) and the
Kodaira-Spencer class of \(X\) is \(v\). Let \(W\subset X\) be the 
pullback of the hyperplane divisor of \(\P^g\).
It follows from \(c_1(W_0)/n = c_1(D)\in H^2(X_0, \BZ)\)
that \(c_1(W_t)/n\in H^2(X_t, \BZ)\) and hence it is a Hodge class.
And since
\begin{equation}\label{lem0:e14}
\Pic(X_t)\isom H^{1,1}(X_t)\cap H^2(X_t, \BZ)
\end{equation}
for K3 surfaces, \(W_t/n\in \Pic(X_t)\). In addition, since
\(\Pic(X_t)\) is torsion free by \eqref{lem0:e14},
\(W_t/n\) is unique in \(\Pic(X_t)\). Hence \(W\sim_{\mathrm{lin}}
nY\) for some divisor \(Y\subset X\), where \(\sim_{\mathrm{lin}}\) is
the linear equivalence. Obviously, \(Y_0\sim_{\mathrm{lin}} D\) and
since \(h^0(\CO_{X_t}(Y_t)) = h^0(\CO_{X_0}(Y_0))\), \(Y\) can be chosen
such that \(Y_0 = D\). We are done.

\begin{rem}\label{apdx1:rem1}
Let \(\m_g\) be the moduli space  of K3 surfaces of genus \(g\).
\lemref{lem0} says that every connected component of \(\m_g\)
is smooth of dimension \(19\). Therefore, we obtain an elementary
proof for this well-known result, which was originally proved
using transcendental methods. See also \cite{CLM} for another elementary
proof of \(\dim \m_g = 19\).

On the other hand, it is also known from the transcendental theory of K3
surfaces that \(\m_g\) is irreducible. Note that the irreducibility of
\(\m_g\) is fundamental to our degeneration argument,
since we rely on the very fact that every K3 surface can be deformed to
a BL K3 surface.
However, it does not seem to be any way of avoiding the use of deep
transcendental theory in order to assert the irreducibility of
\(\m_g\). 
\end{rem}

\section{Recovery of the Counting Formula of Yau-Zaslow-Bryan-Leung}
\label{apdx2}

Let \(N_g\) be the number of rational curves in the primitive class of
a general K3 surface of genus \(g\). We are trying recover the following
remarkable formula of Yau-Zaslow \cite{Y-Z} and Bryan-Leung
\cite{B-L}:
\begin{equation}\label{apdx2:e1}
\sum_{g=0}^\infty N_g q^g = \frac{q}{\Delta(q)} = \prod_{n=1}^\infty
(1-q^n)^{-24}
\end{equation}
where we let \(N_0 = 1\) and \(N_1 = 24\).

By the analysis in \secref{s4}, it is not hard to see the following:

\begin{prop}\label{apdx2:prop1}
Each possible configuration of the stable reduction \(\wt{Y}_0\)
counts exactly one for \(N_g\).
\end{prop}

The above proposition is not hard to prove but it is quite tedious to
write down the whole argument. Basically, by the analysis in
\secref{s4}, \(Y_t\) has exactly \(m\) nodes in the neighborhood of
\(E\) if \(Y_0\) contains \(E\) with multiplicity \(m\); these
\(m\) nodes approach the points \(\nu(\overline{w_{i-1} u_i}) \cap N\)
and \(p_\alpha\) as \(t\to 0\). In order to prove
\propref{apdx2:prop1}, one just has to show that the points
\(\nu(\overline{w_{i-1} u_i}) \cap N\) and \(p_\alpha\) can be
deformed to \(m\) nodes on the general fiber in a ``unique'' way. See
e.g. \cite{CH1}, \cite{CH2}, \cite{CH3} and \cite{C1} for
how to carry out this line of argument. We will leave the details to
the readers.

So it suffices to count the number of possible configurations
of \(\wt{Y}_0\) according to the description given at the end of
\secref{s4}. The number of possible configurations of \(\wt{Y}_0\)
over \(E\) is the same as the number of the sequences \(\{\mu,
\lambda_i\}\) satisfying \eqref{s4:e22}, \eqref{s4:e23} and
\eqref{s4:e24}. We claim that

\begin{prop}\label{apdx2:prop2}
There are exactly \(P(m)\) sequences \(\{\mu,
\lambda_i\}\) satisfying \eqref{s4:e22}, \eqref{s4:e23} and
\eqref{s4:e24}, where \(P(m)\) is the partition number of \(m\), i.e.,
\begin{equation}\label{apdx2:e2}
\prod_{n=1}^\infty (1-q^n)^{-1} = \sum_{m=0}^\infty P(m) q^m.
\end{equation}
\end{prop}

Assume that \propref{apdx2:prop2} holds and then the total number of
possible configurations of \(\wt{Y}_0\) is
\begin{equation}\label{apdx2:e3}
\sum_{m_1+m_2+...+m_{24} = g} P(m_1) P(m_2) ... P(m_{24})
\end{equation}
where \(m_1, m_2, ..., m_{24}\) are the multiplicities of \(Y_0\)
along the \(24\) rational nodal curves \(F_1, F_2, ..., F_{24}\in
|F|\). Obviously, the number given by \eqref{apdx2:e3} is the
coefficient of \(q^g\) in the power series \eqref{apdx2:e1}, i.e.,
\(N_g\). So we are done provided we can prove \propref{apdx2:prop2}.

Let
\begin{equation}\label{apdx2:e4}
G(q, z) = (1+z) \prod_{k=1}^\infty \left((1+q^k z) (1 + q^k
z^{-1})\right).
\end{equation}
We claim that

\begin{lem}\label{apdx2:lem1}
The number of the sequences \(\{\mu,
\lambda_i\}\) satisfying \eqref{s4:e22}, \eqref{s4:e23} and
\eqref{s4:e24} is the same as the coefficient of \(q^m\) in the
power series expansion of \(G(q, z)\).
\end{lem}

\begin{proof}
It follows from the correspondence
\begin{equation}\label{apdx2:e5}
\begin{split}
\{\mu, \lambda_i\} \leftrightarrow & (q^{\lambda_1}z)(q^{\lambda_2}z)
... (q^{\lambda_\mu}z)\\
& \cdot (q^{\lambda_{-1} + 1} z^{-1})
(q^{\lambda_{-2} + 1} z^{-1})
... (q^{\lambda_{-\mu} + 1}z^{-1}).
\end{split}
\end{equation}
\end{proof}

Let us write
\begin{equation}\label{apdx2:e6}
G(q, z) = \sum_{d=-\infty}^\infty C_d z^d
\end{equation}
where \(C_d\in \BC[[q]]\). Then by \lemref{apdx2:lem1},
\propref{apdx2:prop2} holds if and only if
\begin{equation}\label{apdx2:e7}
C_0 = \sum_{m=0}^\infty P(m) q^m = \prod_{n=1}^\infty (1-q^n)^{-1}.
\end{equation}
So it remains to verify \eqref{apdx2:e7}. Our strategy is to first
calculate \(C_{0,n}\) as in
\begin{equation}\label{apdx2:e8}
G_n(q,z) = (1+z) \prod_{k=1}^n \left((1+q^k z) (1 + q^k
z^{-1})\right) = \sum_{d=-\infty}^\infty C_{d,n} z^d.
\end{equation}
and then take the limit \(\lim_{n\to\infty} C_{0,n}\) to obtain
\(C_0\). As long as \(|q|<1\) and \(z\ne 0\), this process makes sense
analytically.

Observe that \(G_n(q, z)\) satisfies the functional equation
\begin{equation}\label{apdx2:e9}
(z + q^n) G_n(q, qz) = (1+q^{n+1} z) G_n(q, z).
\end{equation}
This gives a recursion relation on the coefficients \(C_{d,n}\):
\begin{equation}\label{apdx2:e10}
\begin{split}
&\quad q^{d-1} C_{d-1,n} + q^{n+d} C_{d,n} = C_{d,n} + q^{n+1}
C_{d-1,n}\\
& \Leftrightarrow C_{d-1,n} = \frac{1-q^{n+d}}{q^{d-1}(1-q^{n-d+2})}
C_{d,n}
\end{split}
\end{equation}
for \(-n < d < n+2\). Combining this with \(C_{n+1,n} =
q^{n(n+1)/2}\), we obtain
\begin{equation}\label{apdx2:e11}
C_{0,n} = \frac{(1-q^{2n+1})(1 - q^{2n}) ... (1 -
q^{n+2})}{(1-q)(1-q^2) ... (1-q^n)}.
\end{equation}
Obviously, taking the limit \(C_0 = \lim_{n\to\infty} C_{0,n}\) yields
\eqref{apdx2:e7}. This finishes the proof of \propref{apdx2:prop2} and
hence the recovery of the counting formula \eqref{apdx2:e1}.

\begin{rem}\label{apdx2:rem1}
Here we count the sequences \(\{ \mu, \lambda_i\}\). An alternative
way is to count the sequences \(\{\mu, \mu_j\}\) satisfying
\eqref{s4:e3}, \eqref{s4:e6}, \eqref{s4:e9} and \eqref{s4:e20}. Since
\(\{ \mu, \lambda_i\}\) and \(\{\mu, \mu_j\}\) are ``dual'' to each
other (see \figref{F8}), we may regard this as the {\it dual counting\/} of
what we did above and it should give the same number \(P(m)\).
It turns out that the number of the sequences \(\{\mu, \mu_j\}\) satisfying
\eqref{s4:e3}, \eqref{s4:e6}, \eqref{s4:e9} and \eqref{s4:e20} is
given by the coefficient of \(q^m\) in the expansion of
\begin{equation}\label{apdx2:e12}
\sum_{k=0}^\infty \frac{q^{k^2}}{(1-q)^2(1-q^2)^2...(1-q^k)^2}.
\end{equation}
This leads to the combinatorial identity
\begin{equation}\label{apdx2:e13}
\begin{split}
\prod_{n=1}^\infty (1-q^n)^{-1} &= \sum_{k=0}^\infty
\frac{q^{k^2}}{(1-q)^2(1-q^2)^2...(1-q^k)^2}\\
&= 1 + \frac{q}{(1-q)^2} + \frac{q^4}{(1-q)^2 (1-q^2)^2} + ....
\end{split}
\end{equation}
However, we do not know any direct way to derive
\eqref{apdx2:e13}. We believe that \eqref{apdx2:e13} is known. 
If it is not, it remains an interesting question trying to find a
direct proof for it, a proof without resorting to the correspondence
between \(\{\mu, \lambda_i\}\) and \(\{\mu, \mu_j\}\).
\end{rem}

\begin{rem}\label{apdx2:rem2}
Notice that we did not recover the full formula of Bryan and
Leung. They counted the number of not only rational curves but also
genus \(n\) curves in the primitive class passing through \(n\) general
points.
It is possible to recover their full formula along our line of
argument, but some extra work is needed. The basic setup is the
following. Let \(Y\subset X\) be a family of genus \(n\) curves in the
primitive class of \(X_t\) passing through \(n\) fixed
points in general position. Let \(x_1, x_2,..., x_n\) be the \(n\)
fixed points on \(X_0 = S\) and let \(G_1, G_2,..., G_n\) be the
fibers of \(S\to\P^1\) passing through the points \(x_1,
x_2,...,x_n\), respectively. Then \(Y_0\) is supported along \(G_i\)
and \(F_j\), i.e.,
\begin{equation}\label{apdx2:e14}
Y_0 = \sum_{i=1}^n a_i G_i + \sum_{j=1}^{24} m_j F_j.
\end{equation}
We have analyzed the behavior of \(Y_t\) in the neighborhood of
\(E = F_j\) and classified all possible configurations of
the stable reduction \(\wt{Y}_0\) over \(F_j\). However, we have not
yet done the same for \(Y\) along \(G_i\), which is required for our
counting. On the other hand, this can be carried out along the same
line of argument as we did for \(F_j\). That is, we will repeatedly
blow up \(X\) along \(G = G_i\) until we obtain a nontrivial ruled
surface \(S_\alpha\) over \(G\) on the central fiber. Then we will
analyze the proper transform of \(Y\) under the blowups in much the
same way as we did in \secref{s4}. It will finally come down to the study
of certain curves on \(S_\alpha\). 
Hopefully, we will do this in a future paper.
\end{rem}

\end{document}